\documentclass[a4paper, 12pt]{amsart}
\usepackage{amsmath, latexsym, amssymb, amscd, amsfonts,
enumerate, url, chapterbib, float, graphicx, wrapfig}
\usepackage{mathrsfs}
\usepackage{hyperref}
\newtheorem{theorem}{Theorem}[section] 
\newtheorem{lemma}[theorem]{Lemma}     
\newtheorem{corollary}[theorem]{Corollary}
\newtheorem{proposition}[theorem]{Proposition}
\newtheorem{conjecture}{Conjecture}
\theoremstyle{remark}
\newtheorem*{remark}{Remark}

\newcommand{\meas}{\operatorname{meas}}
\newcommand{\PSL}{\operatorname{PSL}}
\newcommand{\GCD}{\operatorname{GCD}}

\newcommand{\bC}{\mathbb{C}}

\newcommand{\bR}{\mathbb{R}}
\newcommand{\Q}{\mathbb{Q}}
\newcommand{\bQ}{\mathbb{Q}}
\newcommand{\zed}{\mathbb{Z}}

\newcommand{\uH}{\mathbb{H}}

\newcommand{\fd}{\mathfrak{d}}
\newcommand{\fp}{\mathfrak{p}}
\newcommand{\fa}{\mathfrak{a}}
\newcommand{\fb}{\mathfrak{b}}
\newcommand{\fc}{\mathfrak{c}}
\newcommand{\fl}{\mathfrak{h}}

\newcommand{\cC}{{\mathscr{C}}}

\newcommand{\cF}{{\mathscr{F}}}

\newcommand{\cO}{{\mathscr{O}}}
\newcommand{\cM}{{\mathscr{M}}}

\newcommand{\sq}{{\mathrm{sq}}}

\newcommand{\supp}{{\mathrm{supp }}}
\newcommand{\kr}{{\mathrm{kr}}}

\title[Equidistribution of bounded torsion CM points]
 {Equidistribution of bounded torsion CM points} %

\author{B. Hough}
\address{Department of Mathematics, Stanford University, 450 Serra Mall,
Stanford, CA, 94305, USA}
\curraddr{Institute for Advanced Study, Princeton, NJ 08540 USA.}
\email{hough@math.ias.edu}


\thanks{This material is based upon work supported by the National Science
Foundation under agreement No.\ DMS-1128155. Any opinions, findings and
conclusions or recommendations expressed in this material are those of the
author and do not necessarily reflect the views of the National Science
Foundation.}

\thanks{The author  acknowledges the support of a Ric Weiland
Fellowship during his graduate studies at Stanford University.}

\begin{document}
\maketitle

\begin{abstract}
 Averaging over imaginary quadratic fields, we prove, quantitatively,  the
equidistribution of CM points associated to 3-torsion classes in the class
group.  We conjecture that this equidistribution holds for points associated to
ideals of any fixed odd order.  We prove a partial equidistribution result in
this direction and present empirical evidence.
\end{abstract}

\section{Introduction}
 \begin{figure}[h!]
\includegraphics[width=.15 \textwidth]{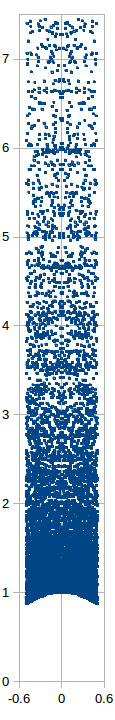}
\caption{Heegner points associated to 101-torsion
classes in imaginary quadratic fields of
discriminant $\approx -4 \cdot 10^6$.}\label{101_figure}
\end{figure}

 \begin{figure}[h!]
\includegraphics[scale = .5]{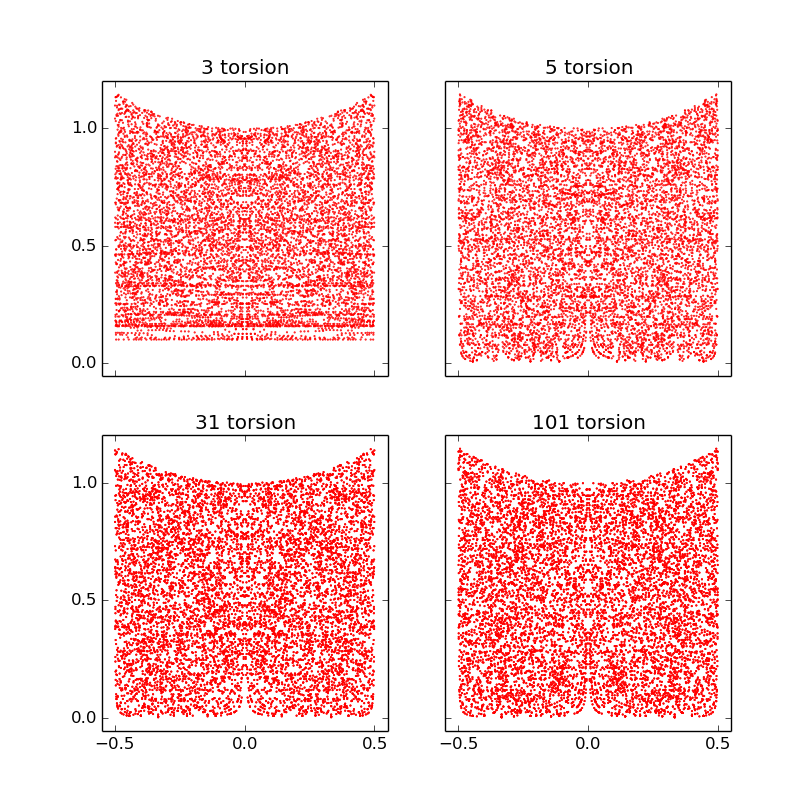}
\caption{Heegner points associated to torsion classes of fixed order
 in imaginary quadratic fields of
discriminant $\approx -4 \cdot 10^6$.  The transformation $y
\mapsto \frac{1}{y}$ has been made, so that the ambient
measure is Lebesgue.}\label{torsion_figure}
\end{figure}
Let $-D$, $D > 0$ be a fundamental discriminant, and write $H_k(-D)$ for the
order $k$ elements in the class group $H(-D)$.  Probably the easiest-to-state
consequence of the  Cohen-Lenstra Heuristics \cite{cohen_lenstra} for imaginary
quadratic fields is the prediction that when fields are ordered by increasing
size of discriminant, for any odd $k > 1$ the average of $|H_k|$ is
asymptotically 1,\footnote{We use $\flat$ to restrict sums to fundamental
discriminants.}
\begin{equation}\label{c-l_prediction}
 {\sum_{ 0 < D < X}}^\flat |H_k(-D)| \sim {\sum_{0 < D < X}}^\flat 1, \qquad X
\to \infty .
\end{equation}
 At any rate, in the special case $k = 3$ this is the only evidence for the
Heuristics which is actually known, thanks to a theorem of Davenport and
Heilbronn \cite{davenport_heilbronn}.
 We wish to broaden the prediction (\ref{c-l_prediction}) to the assertion that
the shapes of lattices of any given odd torsion appear with a common uniform
intensity among the shapes of  all two dimensional lattices, as the discriminant
grows. We will see that this broader interpretation helps to explain the
discrepancy between (\ref{c-l_prediction}) and tabulated data.

 To elaborate, an ideal $\fa$ in the imaginary quadratic field $\bQ(\sqrt{-D})$
is a two-dimensional lattice in $\bC$. To this lattice attach a complex
number $z_\fa$, which is the ratio of any two of it's generators; this number
characterizes the shape of the ideal up to  homothety.  After possibly
exchanging the role of the generators, $z_\fa$ is in the upper half plane
$\uH$, and making a linear change of basis,  it lies
on the modular surface $\cF = \PSL_2(\zed) \backslash \uH$.  This point is
common to all ideals of the same shape (ideal class), and is the CM point of the
class.  Now a famous theorem of Linnik \cite{linnik} and Duke \cite{duke}
asserts that the CM points of classes in $H(-D)$ equidistribute with respect to
the translation-invariant hyperbolic probability measure \[d\mu(z) =
\frac{3}{\pi} \frac{dx dy}{y^2}\] on $\cF$, as $D \to \infty$. Motivated by the
Linnik-Duke theorem, we make our  conjecture.

 \begin{conjecture}\label{equidistribution_conjecture} Let
$K$ be a continuous function of compact support on $\cF$.  For each odd $k > 1$
we have
\[
 \lim_{X \to \infty} {\sum_{0 < D< X}}^\flat \sum_{[\fa] \in H_k(-D)}
K(z_{[\fa]})\bigg/ {\sum_{0 < D< X}}^\flat 1 = \int_{\cF} K(z)d\mu(z).
\]
\end{conjecture}

\noindent The conjecture is well supported by visual evidence, see Figures
\ref{101_figure} and \ref{torsion_figure}.  In fact, notice that the
equidistribution already suggests itself in ranges of discriminants at which
the convergence in the Cohen-Lenstra Heuristics (\ref{c-l_prediction}) is
unconvincing, see Table \ref{torsion_count_table}.
\begin{tiny}
\begin{table}[h!]
 \label{torsion_count_table}
\begin{tabular}{|l|l|l|l|l|l|l|}
\hline
 Range & \# Disc. & 3 & 5 & 7 &11 & 31\\
 \hline
$	-1024	$	&	53	&	28	&
40	&	36	&	10	&	0	\\
$	-2048	$	&	104	&	80	&
64	&	78	&	60	&	0	\\
$	-4096	$	&	206	&	142	&
172	&	162	&	150	&	0	\\
$	-8192	$	&	415	&	316	&
364	&	336	&	240	&	0	\\
$	-16384	$	&	831	&	632	&
752	&	738	&	650	&	270	\\
$	-32768	$	&	1660	&	1338	&
1544	&	1578	&	1330	&	690	\\
$	-65536	$	&	3320	&	2730	&
3192	&	2850	&	2770	&	1890	\\
$	-131072	$	&	6638	&	5532	&
6200	&	6276	&	5800	&	4860	\\
$	-262144	$	&	13286	&	11480	&
12844	&	12348	&	12110	&	10830	\\
$	-524288	$	&	26558	&	23254
&	25072	&	25614	&	25840	&	21210	\\
$	-1048576	$	&	53114	&
47144	&	51328	&	51960	&	50540	&	45210	\\
$	-2097152	$	&	106251	&
95716	&	102340	&	104724	&	103170	&	96960	\\
$	-4194304	$	&	212485	&
193416	&	208288	&	210108	&	207290	&	195570	\\
$	-8388608	$	&	424972	&
391050	&	417516	&	418248	&	415590	&	398550	\\
$	-16777216	$	&	849944	&
789452	&	836176	&	838776	&	832600	&	815790	\\
$	-33554432	$	&	1699872	&
1592438	&	1675940	&	1683882	&	1675150	&	1645380	\\
$	-67108864	$	&	3399779	&
3208270	&	3363532	&	3383604	&	3361140	&	3324120	\\
$	-134217728	$	&	6799584	&
6459970	&	6736896	&	6761478	&	6765140	&	6685350	\\
$	-268435456	$	&	13599079
&	12988450	&	13484300	&	13582980	&
13555960	&	13422870	\\
$	-536870912	$	&	27198220
&	26116790	&	27013804	&	27078228	&
27113010	&	26934720 \\
\hline
\end{tabular}
\caption{Discriminants of the form $4d$, $d \equiv 2 \bmod 4$ are counted in
each specified diadic range, between $2R$ and $R$.  The counts appearing below
each prime $p$ are the
corresponding counts of order $p$ class group elements.}
\end{table}
\end{tiny}

Our main result is a quantitative proof of Conjecture
\ref{equidistribution_conjecture} for the case $k = 3$.
We also have a partial result toward the conjecture for larger $k$, which
asserts that  the CM points are equidistributing `in the cusp'.   It is a
confounding fact that, at least on the basis of visible evidence, the cusp
appears to be the last place where the CM points equidistribute.

\section*{Notation and conventions}

All limiting statements are taken with respect to a growing parameter $X$,
which is a bound for the size of discriminants considered. For positive
functions $A(X)$, $B(X)$, $A \sim B$ means $\lim \frac{A}{B}=1$. We use the
Vinogradov notation $A \ll B$ with the same meaning as $A = O(B)$. $A \asymp B$
means $A \ll B$ and $B \ll A$.
 $\epsilon$ is reserved for a fixed positive parameter which may be taken
arbitrarily small.

Given integrable function $f$ on $\bR^+$, its Mellin transform is defined,
where absolutely convergent, by
\[
 \tilde{f}(s) = \int_0^\infty f(x)x^{s-1}dx, \qquad s \in \bC
\]
and possibly extended elsewhere by analytic continuation.

\section{Precise statement of results}
Recall that the ring of integers in an imaginary quadratic field takes  one of
three forms depending on the behavior  of the discriminant $-D$ at the prime 2.
Since we perform calculations in the ring of integers, for the
remainder of this article we restrict to fundamental discriminants of the form
$-D = -4d$ where $d > 0$, $d \equiv 2 \bmod 4$ is square-free; all of our
arguments carry
over to the other two cases with minor modifications.  In this case, the ring
of
integers is given by $\cO = \zed[\sqrt{-d}]$ within the field $\bQ(\sqrt{-d})$.

We build on the earlier work of Soundararajan \cite{sound}, and much earlier,
Akeny and Chowla \cite{akeny_chowla}, who studied the divisibility problem for
the class group through a parameterization of \emph{primitive ideals}, see also
\cite{byeon_koh} in the real quadratic setting, and \cite{heath-brown} for the
best result on divisibility by 3.  A primitive ideal $\fa \subset \cO$ is an
ideal that does not admit a factorization $\fa = (p \cO) \cdot \fb$ where $p$ is
a prime of $\zed$ and $\fb$ is another ideal of $\cO$.  At the level of
lattices, this says that $\fa$ is not an integer dilation of another ideal.  A
useful characterization of the primitive ideals is that they are exactly those
ideals $\fa$ for which $\{0, 1, ..., N\fa -1\}$ forms a complete set of residues
for $\cO/\fa$.\footnote{For another characterization in terms of the  prime
factorization, see Section \ref{parametrization_section}.}  In particular, this
means that there is a canonical choice of generators for the lattice $\fa$ given
by $\fa = [N\fa, b + \sqrt{-d}]$ where $b$ is uniquely determined by the
conditions \[-\frac{N\fa}{2} < b \leq \frac{N\fa}{2}, \qquad b \equiv
-\sqrt{-d} \bmod \fa.\]
To $\fa$ is then associated the  `Heegner point'
\begin{equation}\label{def_heegner_point}
 z_\fa = \frac{b + \sqrt{-d}}{N\fa}.
\end{equation}
Note that this point lies in the strip $\left(-\frac{1}{2}, \frac{1}{2}\right]
\times \bR^+ = \Gamma_\infty
\backslash \uH$,  where $\Gamma_\infty$ is the subgroup of $\Gamma$ stabilizing
the cusp $\infty$. It is a pretty geometric fact that, fixing an ideal class
$[\fa]$ in the class group $H(-D)$,  the collection of Heegner points of
primitive ideals of class $[\fa]$ are exactly the images of the CM point
$z_{[\fa]}$ in the various fundamental domains for $\Gamma \backslash \uH$
within the strip $\Gamma_\infty \backslash \uH$ (see \cite{iwaniec_kowalski},
Chapter 22).  Therefore, the equidistribution of CM points within the
fundamental domain $\cF$ is equivalent to the equidistribution of the
corresponding Heegner points in the strip $ \Gamma_\infty \backslash \uH$, and
this is the point of view that we shall adopt.  We also introduce the notation
$P_k(-D)$ to denote the primitive ideals with classes in $H_k(-D)$.

Our first result establishes the equidistribution for 3-torsion Heegner points.
\begin{theorem}\label{3_part_equidistribution_theorem}
 Let $K(x,y)$ be a continuous function, compactly supported in the strip
$
\Gamma_\infty \backslash \uH,$ and let $\phi: \bR^+ \to
\bR^+$ be a smooth function of compact support. Let $T = T(X)$ be a parameter
satisfying $1 \leq T \leq
X^{\frac{1}{6}-\epsilon}$.
Then, as $X \to \infty$,
\begin{align*}
  &\sum_{\substack{d \equiv 2 \bmod 4\\ \text{square-free}}}
\phi\left(\frac{d}{X}\right) \sum_{ \fa \in P_3(-4d)  } K\left(\Re z_{\fa},
\frac{\Im
z_{\fa}}{T}\right)\\& \sim
\sum_{\substack{d \equiv 2 \bmod 4\\ \text{square-free}}}
\phi\left(\frac{d}{X}\right)
\int_{\Gamma_\infty \backslash \uH} K\left(x,
\frac{y}{T}\right) \frac{3}{\pi}\frac{dx dy}{y^2}.
\end{align*}
\end{theorem}

Notice that this Theorem  gives more than just the
equidistribution in
$\mathscr{F}$, which follows from the case $T=1$, since it also holds
effectively into the cusp, for
$T < X^{\frac{1}{6}-\epsilon}$.   Actually the result is
stronger, still, since we
have given only a qualitative statement, whereas we can
actually give
quantitative estimates with power saving error terms, see
discussion before
Theorem \ref{quantitative_3_torsion_poincare} in the next section.  For
instance, with  discriminants counted with a smooth weight as above, our method
is strong enough to
yield the Davenport-Heilbronn Theorem ((\ref{c-l_prediction}), $k = 3$) with a
negative secondary main term of size $X^{\frac{5}{6}}$,  giving an alternative
proof of
a recent result of Taniguchi-Thorne \cite{tanaguchi_thorne} and Bhargava,
Shankar and Tsimerman \cite{bhargava_shankar_tsimerman}.
Previously Terr
\cite{terr} has  considered a related equidistribution problem for orders
in cubic fields, by
a different method, but his work yields only the qualitative equidistribution.
Since the completion of this work, Terr's result has been further generalized
by Bhargava and Harron to give an analogous result for the shapes of orders in
quartic and quintic fields \cite{BH13}.

For $k > 3$ we cannot prove the full equidistribution, but we can prove that
Heegner points equidistribute `in the cusp'.

\begin{theorem}\label{k_part_equidistribution}
 Let $K$ and $\phi$ as in the previous theorem, and now assume that $k$ is odd,
$k >3$.  Let
$T = T(X)$ be a parameter growing with $X$ in such a way
that
$
 X^{\frac{1}{2} - \frac{1}{k-2}+ \epsilon} < T <
X^{\frac{1}{2} - \frac{1}{k} -\epsilon}.
$
Then, as $X \to \infty$,
\begin{align*}
 & \sum_{\substack{d \equiv 2 \bmod 4\\ \text{square-free}}}
\phi\left(\frac{d}{X}\right) \sum_{\substack{ \fa\in
P_k(-4d)}} K\left(\Re z_{\fa}, \frac{\Im
z_{\fa}}{T}\right)\\ &\sim \sum_{\substack{d \equiv 2
\bmod 4\\ \text{square-free}}} \phi\left(\frac{d}{X}\right)   \int_{\mathscr{F}}
K\left(x,
\frac{y}{T}\right) \frac{3}{\pi}\frac{dx dy}{y^2}.
\end{align*}
\end{theorem}

\begin{corollary}
 For any odd $k \geq 5$, as $X \to \infty$
 \[
  {\sum_{D < X}}^\flat |H_k(-D)| \gg X^{\frac{1}{2} + \frac{1}{k-2} -\epsilon}.
 \]
\end{corollary}
In the case $k = 5$ this improves the bound $ {\sum_{D < X}}^\flat |H_5(-D)|
\gg X^{\frac{4}{5}}$ from \cite{sound}.

The reader will no doubt have noticed that in both theorems we no longer claim
the the equidistribution of $k$-torsion Heegner points in  the cusp once $\Im(z)
> X^{\frac{1}{2} - \frac{1}{k}}$.  There is a good reason for this -- see Figure
\ref{torsion_figure}.  If $\fa$ is a primitive $k$-torsion ideal in
$\zed[\sqrt{-d}]$ then $\fa^k = (x + y \sqrt{-d})$ is principal, and $y \neq 0$,
since $\fa$ is primitive.  Hence $N\fa^k = x^2 + dy^2 \geq d$ so that we have
the upper bound
\[
 \Im(z_{\fa}) = \frac{\sqrt{d}}{N\fa} \leq d^{\frac{1}{2} - \frac{1}{k}}.
\]
Since the set $\left\{z \in \Gamma_\infty \backslash \uH: \Im(z) >
X^{\frac{1}{2} - \frac{1}{k}}\right\}$ has hyperbolic volume $\asymp
X^{-\frac{1}{2} + \frac{1}{k}}$,  the absence of Heegner points in this set
suggests a negative secondary term in (\ref{c-l_prediction}) of size
$X^{\frac{1}{2} + \frac{1}{k}}$.  After a fashion, we are able to determine this
quantity of missing torsion points as the negative secondary main term in the
following theorem.

\begin{theorem}\label{k_torsion_secondary_term}
 Let $k> 3$ be odd.  Let $\phi, \psi$ be $C^\infty$ functions on
$\bR^+$ with $\phi$ having compact support, and $\psi$
supported in $[1, \infty)$, with $\psi \equiv 1$ on a
neighborhood of $\infty$. Denote $\tilde{\phi}, \tilde{\psi}$ the Mellin
transforms.
There exists a $\delta = \delta_k > 0$ such that for $T$ in the range
\[X^{\frac{1}{2}-\frac{1}{k} - \delta_k } \leq T \ll X^{\frac{1}{2} -
\frac{1}{k}}\]
we have the asymptotic with two main terms
\begin{align*}
 &\sum_{\substack{d \equiv 2\bmod 4\\ \text{square-free}}}
\phi\left(\frac{d}{X}\right) \sum_{\substack{\fa \in
P_k(-4d)}}\psi\left(\frac{\Im
z_{[\fa]}}{T}\right) \\&= \frac{6}{\pi^3}
\tilde{\phi}(1)\tilde{\psi}(-1)\frac{X}{T} + c_{k}
\tilde{\phi}\left(\frac{1}{2} + \frac{1}{k}\right)
X^{\frac{1}{2} + \frac{1}{k}} + o\left(X^{\frac{1}{2} +
\frac{1}{k}}\right);\\
c_{k} &=
\frac{\Gamma(\frac{1}{2}- \frac{1}{k})\zeta\left(1-\frac{2}{k}\right)}{k
\pi^{\frac{3}{2}}\Gamma(1
-
\frac{1}{k})} \\& \times \left[1 -
2^{\frac{1}{k}} + 2^{1-\frac{1}{k}}\right]\prod_{p \text{ odd}}
\left[1 + \frac{1}{p+1}\left(\frac{1}{p^{\frac{1}{k}}} -
\frac{1}{p^{1-
\frac{2}{k}}} - \frac{1}{p^{1-\frac{1}{k}}}-
\frac{1}{p}\right)\right].
\end{align*}
The secondary term of size $X^{\frac{1}{2} + \frac{1}{k}}$ is negative, since
$\zeta\left(1 - \frac{2}{k}\right) < 0$.
\end{theorem}
\begin{remark}
 Our proof will show that we may take any $\delta_k < \frac{2}{k^2}$.
\end{remark}
When $k = 3$, the term $c_3 \tilde{\phi}(\frac{5}{6}) X^{\frac{5}{6}}$  is the
actual negative secondary term in the Davenport-Heilbronn Theorem when
discriminants are counted with smooth weight $\phi$.  For $k = 5, 7$, inclusion
of this secondary term in the right side of (\ref{c-l_prediction}) brings this
prediction into good agreement with tabulated data for relatively small
discriminants, see Table \ref{57_torsion_table}.  For $k \geq 9$, the
agreement is not as good in the region in which we have numerical data.

\begin{tiny}
\begin{table}[h!]
 \label{57_torsion_table}
\begin{tabular}{|l|l|l|l|l|l|}
\hline
X & $\sum_{d<X} h_5(-4d)$ & Cohen-Lenstra & CL-$\Sigma$ & CL + $c_5
X^{\frac{7}{10}}$ & CL + $c_5
X^{\frac{7}{10}} - \Sigma$\\
\hline
1000000	&	194464	&	202642	&	8178	&	194510	&
46	\\
2000000	&	392996	&	405285	&	12289	&	392074	&
-922	\\
4000000	&	791328	&	810569	&	19241	&	789108	&
-2220	\\
8000000	&	1588520	&	1621139	&	32619	&	1586275	&
-2245	\\
16000000	&	3186224	&	3242278	&	56054	&	3185641
&	-583	\\
32000000	&	6393960	&	6484556	&	90596	&	6392548
&	-1412	\\
64000000	&	12818136	&	12969112	&	150976
&	12819645	&	1509	\\
128000000	&	25673816	&	25938223	&	264407
&	25695414	&	21598	\\
\hline
\multicolumn{6}{l}{} \\
\hline
X & $\sum_{d<X} h_7(-4d)$ & Cohen-Lenstra & CL-$\Sigma$ & CL + $c_7
X^{\frac{9}{14}}$ & CL + $c_7
X^{\frac{9}{14}} - \Sigma$\\

\hline
1000000	&	197094	&	202642	&	5548	&	196900	&
-194	\\
2000000	&	397902	&	405285	&	7383	&	396318	&
-1584	\\
4000000	&	796266	&	810569	&	14303	&	796568	&
302	\\
8000000	&	1595088	&	1621139	&	26051	&	1599277	&
4189	\\
16000000	&	3201048	&	3242278	&	41230	&	3208143
&	7095	\\
32000000	&	6427098	&	6484556	&	57458	&	6431257
&	4159	\\
64000000	&	12870768	&	12969112	&	98344
&	12885890	&	15122	\\
128000000	&	25832964	&	25938223	&	105259
&	25808279	&	-24685	\\
\hline
\end{tabular}
\caption{Aggregate order 5 and 7 elements in the class group of quadratic fields
of discriminant
$-4d$, $d < X$ are tabulated.  Conjectural secondary main terms of size
$X^{\frac{7}{10}}$ and $X^{\frac{9}{14}}$ respectively improve the numerical
fit of the Cohen-Lenstra heuristics.}
\end{table}
\end{tiny}

\section{Discussion of method}\label{method_section}
One description of the divisibility argument in \cite{sound} is that the norm
equation
\[
 N \fa^k = m^k = x^2 + dy^2
\]
is used to parametrize and count some $k$-torsion primitive ideals of
$\zed[\sqrt{-d}]$ within a band  in the cusp of $\cF$.  We refine the
parameterization used so as to give the exact location of the counted points.  A
precise statement of the parameterization along with a local version is at the
beginning of the
next section.

Our proofs of equidistribution are by Weyl's criterion, that is, we use that
the
linear span of functions of the form
\[
  e(f x)\psi(y), \qquad f \in \zed,\; \psi \in C_c^\infty(\bR^+)
\]
is dense in the space of continuous functions of compact support on the strip
$\bR/\zed \times \bR^+$. This reduces the proofs of Theorems
\ref{3_part_equidistribution_theorem} and \ref{k_part_equidistribution} to the
estimates (here $\psi_T(y) = \psi(\frac{y}{T})$, and $\tilde{\phi}$ and
$\tilde{\psi}$ denote the Mellin transforms)
\begin{align}\label{weyl_criterion}
& \sum_{\substack{d \equiv 2 \bmod 4\\ \text{square-free}}}
\phi\left(\frac{d}{X}\right)  \sum_{\fa \in P_k(-4d)} e(f \Re z_\fa) \psi_T(\Im
z_\fa)\\&\notag = \delta_{f=0} \tilde{\phi}(1) \tilde{\psi}(-1)\frac{X}{T}  +
o\left(\frac{X}{T}\right), \end{align}
for any $\phi, \psi \in C_c^\infty(\bR^+)$, $f \in \zed$ and for $T$ in the
stated ranges of the theorems.
Strictly speaking, to obtain quantitative equidistribution
 one requires estimates of  the type (\ref{weyl_criterion}) with error terms
that make explicit the dependence on the frequency $f$ and function space norms
of
the test function $\psi$.  In the quantitative theorems that we state below we
have tracked the frequency dependence but omit the dependence on $\psi$.

 \begin{theorem}\label{quantitative_3_torsion_poincare}
 Let $\phi, \psi \in C_c^\infty(\bR^+)$, with $\phi$ the
function of Theorem \ref{3_part_equidistribution_theorem}.
Let $f \in \zed$ and  $T = T(X)$ be a
parameter that satisfies $1 \leq T \leq X^{\frac{1}{6}-\epsilon}$.
Define $\psi_T(y) = \psi \left(\frac{y}{T}\right)$.  We have
\begin{align*}
 &\sum_{\substack{d \equiv 2 \bmod 4\\ \text{square-free}}}
\phi\left(\frac{d}{X}\right) \sum_{ \substack{\fa \in
P_3(-4d)\\
}} e(f \Re z_\fa) \psi_T ( \Im z_\fa)\bigg/\sum_{\substack{d
\equiv 2 \bmod 4\\
\text{square-free}}} \phi\left(\frac{d}{X}\right) \\ &=
\delta_{f = 0} \cdot \left[
\frac{3}{\pi T} \int_0^\infty \psi(y) \frac{dy}{y^2}  \right] +
O\left((1 + |f|)^{\frac{1}{2}}
\frac{X^{-\frac{1}{8}+\epsilon}}{T^{\frac{5}{4}}}\right) +
O\left(X^{-\frac{1}{6} + \epsilon}\right).
\end{align*}
\end{theorem}
\noindent To obtain Theorem
\ref{3_part_equidistribution_theorem}, approximate the
function $K(x,y)$ as a linear combination of functions
$\psi(y)e(fx)$ and apply the above theorem term-by-term.

For $k$-torsion with $k > 3$ the
estimate that we prove is as follows.
\begin{theorem}\label{k_torsion_poincare}
Let $k > 3$ odd,  $\phi, \psi \in C^\infty(\bR^+)$ with
$\phi$ of compact
support and $\psi$ supported in $[1, \infty)$ with $\psi
\equiv 1$ on a
neighborhood of $\infty$.
Let $f \in \zed$ and let $T = T(X)$
be a parameter, with
$\psi_T(y) = \psi\left(\frac{y}{T}\right)$ as before.  If $f
= 0$ then for $T$
in the range $X^{\frac{1}{2} - \frac{1}{k-2} +\epsilon} < T < X^{\frac{1}{2}
-\frac{1}{k}-\epsilon}$ we have
the asymptotic
\begin{align}\notag
 \sum_{\substack{d \equiv 2\bmod 4\\\text{square-free}}}
&\phi\left(\frac{d}{X}\right) \sum_{\fa \in P_k(-4d)}\psi_T(\Im z_\fa)
\bigg/\sum_{\substack{d \equiv 2
\bmod
4\\ \text{square-free}}} \phi\left(\frac{d}{X}\right)\\=&
\frac{3}{\pi T}
\int_0^\infty \psi(y)\frac{dy}{y^2} + \frac{\pi^2}{2} c_{k}
\frac{\tilde{\phi}\left(\frac{1}{2} +
\frac{1}{k}\right)}{\tilde{\phi}(1)}
X^{\frac{1}{k}-\frac{1}{2}  }   \label{k_torsion_asymptotic}\\&\notag +
 O\left(\frac{X^{\frac{k}{4}-1+
\epsilon}}{T^{\frac{k}{2}}}\right) +
O\left(X^{ \frac{1}{2k-2}-\frac{1}{2}+\epsilon}\right).
\end{align}
with $c_k$ the constant of Theorem
\ref{k_torsion_secondary_term}.

If $f \neq 0$ then for $T$ in the range $X^{\frac{1}{2} -
\frac{1}{k-2} + \epsilon} < T < X^{\frac{1}{2} -\frac{1}{k} -\epsilon}$ we have
the bound
\begin{align*}
 \sum_{\substack{d \equiv 2\bmod 4\\\text{square-free}}}
\phi\left(\frac{d}{X}\right) &\sum_{\substack{\fa \in P_k(-4d) }}e(f \Re z_\fa)
\psi_T(\Im z_\fa) \bigg/\sum_{\substack{d \equiv 2
\bmod
4\\ \text{square-free}}} \phi\left(\frac{d}{X}\right) \\& =
O\left(\frac{X^{\frac{k}{4}-1+\epsilon}}{T^{\frac{k}{2}}}
\right) +
O\left(\frac{|f|^{\frac{1}{2}}X^{\frac{k}{8} -
\frac{1}{2}+\epsilon}}{T^{\frac{k}{4} + \frac{1}{2}}}\right)
+
O\left(X^{ \frac{1}{k}-\frac{1}{2} + \epsilon}\right).
\end{align*}
\end{theorem}
In the range $T > X^{\frac{1}{2} - \frac{1}{k} -\frac{2}{k^2}+\epsilon}$, the
expression
(\ref{k_torsion_asymptotic}) is an asymptotic formula with
two main terms, and
so we obtain Theorem \ref{k_torsion_secondary_term}.  Notice
that the
terms with fixed $f \neq 0$ are dominated by the main term with
$f = 0$ once
$T > X^{\frac{1}{2} - \frac{1}{k-2}+\epsilon}$.
Although we
have stated this Theorem for $\psi$ with $\lim_{t \to
\infty} \psi(t) = 1$, any
function $\psi_0$ with compact support on $\bR^+$ is the
difference of two such
functions.  Thus we may obtain the stated result for any
$\psi$ having compact
support, but there will be no secondary main term. In
particular, Theorem
\ref{k_part_equidistribution} follows from this Theorem by
approximating
$K(x,y)$ in the space of functions of form $\psi(y)e(fx)$.

The remainder of the  paper is concerned with proving
Theorems
\ref{quantitative_3_torsion_poincare} and
\ref{k_torsion_poincare}.

\section{Parameterization}\label{parametrization_section}
The starting point is the following
parameterization of ideals in $k$-torsion classes  of the
class group of $ \bQ(\sqrt{-d})$.

\begin{proposition}\label{global_parametrization} Let $d
\equiv 2 \bmod 4$ be square-free
  and $k \geq 3$ be odd.  The set
$$\{(\ell,m,n,t)\in (\zed^+)^4: \ell m^k = \ell^2n^2 +
t^2d, (\ell mn,t) = 1\}$$ is in bijection with primitive
ideal pairs
$\{\fa, \overline{\fa}\}$
with $\fa \neq 1$ and $\fa^k$ principal in $\Q(\sqrt{-d})$.
Explicitly, the
ideals $\fa, \overline{\fa}$ are
given as $\zed$-modules by $$\fa = [\ell m, \ell nt^{-1} +
\sqrt{-d}], \qquad \overline{\fa} = [\ell m, -\ell nt^{-1} +
\sqrt{-d}]$$ where $N\fa =
\ell m$ and $t^{-1}$ is the inverse of $t$ modulo $m$.  In particular,
\[
 z_\fa = \frac{nt^{-1}}{m} + i \frac{\sqrt{d}}{\ell m}, \qquad
z_{\overline{\fa}} = \frac{-nt^{-1}}{m} + i \frac{\sqrt{d}}{\ell m}.
\]

\end{proposition}

In our statement of results we have already mentioned two characterizations of
the primitive ideals of $\cO$, but for the proof of Proposition
\ref{global_parametrization} it is  convenient to have a third.  Recall
that ideals of $\cO$  have
unique
factorization, with the  behavior in $\cO$ of the primes
$p\mathscr{O}$ of $\zed$ described by the
quadratic
character\footnote{$\left(\frac{\cdot}{p}\right)$ is
  the Legendre symbol.} of $-d$ mod $p$
$$ p \mathscr{O} = \left\{\begin{array}{lll} \mathfrak{p}^2
& & p|d \\
    \mathfrak{p} \overline{\mathfrak{p}} &&
\left(\frac{-d}{p}\right)
    = 1 \\ p \mathscr{O} && \left(\frac{-d}{p}\right) =
    -1\end{array}\right. .$$  We say that $p$ either
ramifies, splits,
or remains inert.  The different is the product of primes
containing
$d$,
\[\mathfrak{d} = \prod_{\mathfrak{p} | (d)} \mathfrak{p}.\]
In this description, an ideal $\mathfrak{a}$ of $\mathscr{O}$ is primitive if
and
only if
it factors as $\mathfrak{a} = \mathfrak{l}\mathfrak{b}$ with
$\mathfrak{l}|\mathfrak{d}$, $(\mathfrak{b},\mathfrak{d}) =
(1)$ and
$(\mathfrak{b}, \overline{\mathfrak{b}}) = (1)$.  In
particular, $\fb$
contains only primes $\fp$ dividing split primes, with at
most one of $\fp,
\overline{\fp}$ appearing.

\begin{proof}[Proof of Proposition \ref{global_parametrization}]
Take $\fa \neq (1)$ primitive with $\fa^k$ principal and
write
$\fa = \fl \fb$ where  $\fl | \fd$ and $(\fb, \fd) = (1)$.
We have $\fb \neq
(1)$ since otherwise $\fa = \fl\Rightarrow [\fl]^k = [\fl] =
[1]$ which forces $\fl = (1)$.
Now
\begin{equation}\label{ideal_factorization}\fa^k\fl^{-(k-1)}
= (x +
t\sqrt{-d})\end{equation} is principal.  It is also
primitive since $(x +
t\sqrt{-d}) = \fl \fb^k$ and $(\fb,\overline{\fb}) = (1)$,
$(\fb,
\fd)=1$.   Let $m = N\fb$, $\ell
=N\fl$ and take norms in eqn. (\ref{ideal_factorization}) to
obtain $\ell m^k = x^2
+ t^2d$.  Here $\ell |x$ so writing $x = \ell n$,  $m^k =
\ell  n^2 + t^2 \overline{\ell}$
where $\ell \overline{\ell} = d$.  Now primitivity of the
ideal $(\ell n + t\sqrt{-d})$
implies $(t, \ell n) = 1$. Also $(m,t) = 1$, since if
$p|(m,t)$ then $p^2|\ell n^2$ so $p|(n, t)$ which is false.
Finally, primitivity of $(\ell n + t\sqrt{-d})$ implies $n,
t \neq
0$.  We may fix $t >0$ by multiplying by $\pm 1$;  the
choice of sign for
$n$ is determined by a choice between the ideals $\fa$ and
$\overline{\fa}$.

Now suppose we begin with a solution $(\ell ,m,n, t)$ to
$\ell m^k = \ell^2n^2 + t^2 d$
with $(\ell mn, t)=1$ and $\ell ,m,n,t>0$.  Observe that
$\ell|d$, so $\ell$ is square-free.  We claim that also
$(m,n)=1$, which implies $(m,d)=1$.  Indeed, $(m,n)=1$
follows from the fact that $d$ is square-free, since if
$p|(m,n)$ then $p\nmid t$ so that $p^2 | d= \frac{\ell m^k -
\ell^2n^2}{t^2}$, a contradiction.

 Write $(\ell n + t\sqrt{-d}) = \fl \fc$
where $\fl|\fd$ and $(\fc, \fd) = 1$.  Then $(\ell)(m^k) =
\fl^2 \fc
\overline{\fc}$
and $(m,d) = 1$ implies $\fl^2 = (\ell)$ and $ \fc
\overline{\fc}=(m^k)$.
Moreover, $\fc$ is primitive since it divides $(\ell n +
t\sqrt{-d})$, and $\fc$
is prime to
$\fd$ so $(\fc, \overline{\fc}) = 1$, and hence there exists
$\fb$ with $\fc =
\fb^k$, $\overline{\fc} = \overline{\fb}^k$.  Note that
$(\fb, \fd) = 1$ and
$\fb$ is primitive.  Then letting $\fa = \fl \fb$,
$\overline{\fa} = \fl
\overline{\fb}$ we get that $\{\fa, \overline{\fa}\} \neq
\{(1),(1)\}$ is a pair of
primitive ideals satisfying $\fa^k =
(\ell)^{\frac{k-1}{2}}(\ell n + t\sqrt{-d})$ is principal.
Since
there were no choices in determining the pair $(\fa,
\overline{\fa})$, this
completes the bijection.

Taking $\fa$ to be the ideal in the pair $(\fa,
\overline{\fa})$ that
corresponds to $n,t>0$, we now specify $\fa$ in terms of
$\ell,m,n,t$. Since
$\fa$ is primitive, $\fa = [N\fa, b + \sqrt{-d}]$ as a
$\zed$-module, where $b$
is
determined modulo $N\fa$.  From the above bijection, $N\fa =
\ell m$, so it remains
to determine $b \bmod{\ell m}$.  Writing $\fa = [\ell m, b +
\sqrt{-d}]$
and multiplying, $$\fa^2 = (\ell )\fb^2 = [\ell^2m^2, \ell
mb + \ell m
\sqrt{-d},
b^2 - d + 2b\sqrt{-d}].$$  For the right side to be
divisible by $\ell$, we must
have $\ell|b^2 - d$ so $\ell|b^2 \Rightarrow \ell |b$ so
write $b = \ell b'$.
Since $\fa$ contains the element $\ell m$, and $\fb^2$
contains both
the elements $\ell m^2$ and $\ell m b' + m \sqrt{-d}$ the
ideal $$\fa
(\fb^2)^{\frac{k-3}{2}} \fb^2 = (\ell)^{-\frac{k-1}{2}}\fa^k
= (\ell n
+ t \sqrt{-d})$$ contains the element $(\ell m)(\ell
m^2)^{\frac{k-3}{2}}(\ell mb'
+ m \sqrt{-d})$.  Hence for some integers $x, y$,
\begin{align*}\ell^{\frac{k+1}{2}}m^{k-1}b' +
\ell^{\frac{k-1}{2}}m^{k-1} \sqrt{-d} &= (\ell n +
t\sqrt{-d})(x + y\sqrt{-d}) \end{align*}
and  therefore
\begin{align*} \ell^{\frac{k+1}{2}}m^k b' +
\ell^{\frac{k-1}{2}}m^k \sqrt{-d} = (\ell n + t\sqrt{-d})(mx
+ my
\sqrt{-d}).\end{align*}  Now factor $\ell m^k = (\ell n +
t\sqrt{-d})(\ell n -
t\sqrt{-d})$ and cancel $(\ell n + t\sqrt{-d})$ from both
sides of the
above equation to find
\begin{align*}&\left(\ell n -
t\sqrt{-d}\right)\left(\ell^{\frac{k-1}{2}}b' +
  \ell^{\frac{k-3}{2}}\sqrt{-d}\right)\\
&\qquad=\left(\ell^{\frac{k+1}{2}}nb'
  +\ell^{\frac{k-3}{2}}td\right)
+\left(\ell^{\frac{k-1}{2}}n-\ell^{\frac{k-1}{2}}
tb'\right)\sqrt{-d} = mx +
my\sqrt{-d}.\end{align*}  Hence $$\ell^{\frac{k-1}{2}}n
\equiv
\ell^{\frac{k-1}{2}}tb'\bmod m \qquad \Rightarrow \qquad b'
\equiv t^{-1}n \bmod m$$ and $b = \ell b'
\equiv \ell nt^{-1} \bmod \ell m$ as claimed. \end{proof}

The above parameterization suggests a local relation of type
\[ m^k = n^2 + t^2 d \qquad \Rightarrow \qquad m^k \equiv
n^2 \bmod t^2.\]  We now give a local parameterization of solutions to this
congruence.
\begin{proposition}\label{local_parametrization}
Let $N>0$ be an integer and $k \geq 1$ be odd.  Define
\[ S_N = \{(m,n) \in ((\zed/N\zed)^\times)^2: m^k\equiv n^2
\bmod
N\}\] and
\[S'_{N}  = \{(m,n) \in ((\zed/4N\zed)^\times)^2 : m^k - n^2
\equiv 2N \bmod
4N\}.\]
The sets $S_N$ and $S_N'$ have the local parameterization \begin{align*} S_N &=
\{(w^2, w^k): w \in (\zed/N\zed)^\times\},\\
S'_{N} &= \{(m+2N, n): (m,n) \in S_{4N}\}.
\end{align*}

Furthermore, given $(m,n)\in \zed^2$, $(mn, N) = 1$ solving $m^k
\equiv n^2
\bmod N^2$, one has the parameterization
\begin{align*}
 &\{(m', n') \in (\zed/N\zed)^2: (m', n') \equiv (m,n) \bmod N, {m'}^k
\equiv {n'}^2 \bmod N^2\} \\&= \{(m + aN, n + a'N): a, a' \in
\zed/N\zed, k a m^{k-1} \equiv 2 a'n \bmod N\}.
\end{align*}
\end{proposition}

\begin{proof}
To prove the parameterization, note that $w \mapsto (w^2,
w^k)$ and
$(m,n)\mapsto m^{-\frac{k-1}{2}}n$ are inverse maps between
$(\zed/N\zed)^\times$
and $S_N$.  The remaining claims are simple modular
arithmetic.
\end{proof}

Ultimately we will solve for $d$ in the parameterization of Proposition
\ref{global_parametrization}, and sieve for $d$ that are fundamental
discriminants.  In bounding the error from the sieve in Section
\ref{sieve_section} we require the following estimate for the number of
primitive ideals of bounded norm in a given ideal class.

\begin{proposition}\label{norm_bound}
Fix an ideal class $[\fa] \in H(-4d)$. Let $Y>0$. We have the
bound
$$\left|\left\{\fb \text{ primitive}: [\fb] = [\fa], N\fb \leq
Y\sqrt{d}\right\}\right| \ll 1 + Y.$$
\end{proposition}

\begin{proof}
The condition $N\fb \leq Y \sqrt{d}$ is equivalent to $\Im
z_\fb  \geq
Y^{-1}$.  Since $z_\fb = \gamma\cdot z_{[\fa]}$ for some
$\gamma \in
\Gamma_\infty \backslash \Gamma$ the result is a consequence
of the simple geometric estimate, valid for any $z$ in the strip $\Gamma_\infty
\backslash \uH$,
\[\left|\left\{\gamma \in \Gamma_\infty \backslash \Gamma:
\Im \gamma z \geq Y^{-1}\right\}\right| \ll 1 + Y,\] see
\cite{iwaniec_spectral} Lemma 2.11.
\end{proof}

We close this section with a bound for certain complete  exponential sums.
Let
\begin{equation}\label{complete_sum} S_k(A, B; q) = {\sum_{w \bmod q}}^\times
e\left(\frac{A
w^{2} + B w^{k}}{q}\right).\end{equation}
This sum
factors as a product
over prime power sums,
\begin{align*}
 &S_k(A, B; q) = \prod_{p^j\| q} S(A \overline{q}_p, B \overline{q}_p
; p^j)\\
 &q_p = \frac{q}{p^j}, \qquad \overline{q}_p q_p \equiv 1 \bmod p^j.
\end{align*}
For the prime power sums we record the following
lemma.
\begin{lemma}\label{S_f_u_lemma}
We have the following evaluation and bounds for $S_k(A, B;
p^n)$.

\begin{enumerate}
  \item[i.] If $p^n |(A,B)$ then $S_k(A, B; p^n) = (p-1)p^{n-1}$.
 \item[ii.] If $p^j \| (A,B)$ with $j < n$ then $S_k(A, B; p^n) =
p^j S_k(\frac{A}{p^j}, \frac{B}{p^j}; p^{n-j})$
  \item[iii.] If $p \nmid (A,B)$ then $|S_k(A,B; p^n)|  \ll_k
p^{\frac{n}{2}}$.
\end{enumerate}
In particular,
\[
 |S_k(A, B; p^n)| \ll_k \GCD(A, B, p^n)^{\frac{1}{2}}p^{\frac{n}{2}}.
\]

\end{lemma}
\begin{proof}
 Items (i) and (ii) are obvious.  
In
case (iii) the bound
holds for $n=1$ by Weil's bounds. For $n > 1$ this is elementary.
\end{proof}

\section{Function notation and properties}\label{function_property_section}
We adopt the following notation regarding Fourier transforms.
For a smooth integrable function $f$ in several variables
denote by
\begin{align*} f^1(u,y, z) & = \int_{-\infty}^\infty
f(x,y,z)e(-ux) dx, \\
 f^2(x,v,z)  &= \int_{-\infty}^\infty f(x,y,z)e(-vy) dy,\\
f^{1,2}(u,v,z) & = \int_{-\infty}^\infty\int_{-\infty}^\infty f(x,y,z)e(-ux -
vy) dx
dy
\end{align*}
the function with the Fourier transform taken
in the first, second,
or both first and second slots.

\begin{lemma}\label{fourier_shift}
Let $F \in \mathscr{S}\left(\bR^2\right)$ be a Schwarz-class
function and set
\[f\left(x,y\right) = F\left(A + Bx, C + Dx + Ey\right),
\qquad B, E \neq
0.\]
Then
\[f^{1,2}\left(u,v\right) = \frac{1}{BE} e\left(\frac{A}{B}u
+ \left(\frac{C}{E}
-
\frac{AD}{BE}\right)v\right) F^{1,2}\left(\frac{1}{B} u -
\frac{D}{BE}v,
\frac{1}{E}v\right).\]
\end{lemma}

Throughout, $\phi \in C_c^\infty(\bR^+)$ is the smooth
function of the
Theorems,  and also fix once for all time a
function $\sigma \in
C_c^\infty(\bR)$ satisfying
\begin{equation}\label{sigma_properties}\sigma \geq
0,\qquad
\supp\left(\sigma\right) \subset [-2,2], \qquad
\sum_{n \in \zed} \sigma \left(n + x\right) = 1.
\end{equation}  Letting
$\sigma^\times(x) =
\sigma(\log x)$ we obtain a related non-negative
function of compact
support on $\bR^+$ satisfying
\begin{equation}\label{sigma_cross_properties}\sum_{n \in \zed}
\sigma^\times\left(e^n x\right)
= 1.\end{equation}

In addition to the fixed $\phi$, let $\{\psi_j\}_{j \in
\zed}$ be
 smooth functions on
$\bR^+$ satisfying uniform support and $C^k$ bounds
 \begin{align}
\label{support_condition}
\supp\left(\psi_j\right) \subset
[e^{-6}, e^6], \qquad \forall j \in \zed, \forall k \geq 0, \qquad
\|\psi_j\|_{C^k} = O_k(1).
\end{align}
 Note in particular that $\phi$ and  $\psi_j$ have
Mellin transforms $\tilde{\phi}, \tilde{\psi}_j$ that are
entire, and satisfy uniform bounds
\begin{equation}\label{mellin_bounds}
 \forall A > 0, \forall s \in \bC^{\times}, \qquad
\left|\tilde{\phi}(s)\right|,\left|\tilde{\psi}_j(s)\right| \ll_A |s|^{-A}.
\end{equation}

For positive parameters $X, Y$, a frequency $f \in \zed$
and a smooth bounded function $\psi$ on $\bR^+$, supported
away from 0, define
\begin{equation}\label{def_Phi_X_Y_a_f}\Phi_{X, Y,f}\left(x,y, z| \psi\right)
=
\phi\left(\frac{x^k -
y^2}{Xz^2}\right)\psi\left(\frac{x^k
  -y^2}{Y^2  x^2z^2}\right)e\left(-\frac{f
y}{xz}\right),\end{equation}
with the interpretation that $\phi$ and $\psi$ vanish at negative argument.
This is the typical function packaging the `Archimedean' data of our analysis.
Also set for $M > 0$ and $F \in \bR$,
\begin{equation}\label{def_Psi_M_F}\Psi_{M,F}\left(x,y|
\psi\right) =
\phi\left(x^k -
y^2\right) \psi\left(\left(x^k -
y^2\right)\frac{M}{x^2}\right)e\left(-\frac{Fy}{x}
\right).\end{equation}  The appropriate $\psi$ will
generally be clear from
the context, in which case the last argument is dropped.  Note
that for fixed
$x$, and for $\psi_j$ satisfying
support condition
(\ref{support_condition}), $\Psi_{M,F}(\psi_j)$ is supported
on
$x^k - y^2
\asymp 1$ so that
\[\meas(\{y: \Psi_{M,F}(x,y|\psi_j) \neq 0 \}) \ll x^{-\frac{k}{2}}.\]  Also,
 $\Psi_{M, F}(\psi_j)$ is
supported on $x \asymp \sqrt{M}$.
In particular,
\begin{equation}\label{L_1_bound}
 \|\Psi_{M,F}(\psi_j)\|_1 \ll M^{-\frac{k-2}{4}}.
\end{equation}

\begin{lemma}The Fourier transforms of $\Phi$ and $\Psi$ are related as follows,
\begin{align}\label{Phi_Psi_FT}\Phi_{X,Y,f}^{1,2}
\left(u,v, z\right) &=
\left(z^2X\right)^{\frac{1}{2} +
\frac{1}{k}}\Psi_{M,F}^{1,2}\left(z^{\frac{2}{k}}X^{
\frac{1}{k}}u,
  zX^{\frac{1}{2}}v\right),\\\notag M
&=\frac{X^{1-\frac{2}{k}}}{Y^2z^{\frac{4}{k}}}, \qquad F = \frac{f
X^{\frac{1}{2}
-
\frac{1}{k}}}{z^{\frac{2}{k}}}.\end{align}
\end{lemma}
\begin{proof}
 We have
\begin{align}\label{Phi_FT}
&\Phi_{X,Y,f}^{1,2}\left(u,v, z\right)\\ \notag &= \int_{\bR^2}
\phi\left(\frac{x^k-y^2}{
  Xz^2}\right)\psi\left(\frac{x^k-y^2}{Y^2x^2z^2}\right)
e\left(-\frac{fy}{xz}-ux - vy\right)dxdy
\\\notag &= \left( Xz^2\right)^{\frac{1}{2} + \frac{1}{k}}\int_{\bR^2}
\phi\left(x^k -
y^2\right)\psi\left(\frac{x^k-y^2}{x^2}\frac{X^{1-\frac{2}{k}}}{
Y^2z^{\frac{4}{k}}}
\right)\\ \notag &\qquad  \times e\left(-\frac{fX^{\frac{1}{2} -
\frac{1}{k}}y}{xz^{\frac{2}{k}}}-\left(Xz^2\right)^{\frac{1}{k}}ux
- \left(Xz^2\right)^{\frac{1}{2}}vy\right)dx dy
\\ \notag&=
\left(Xz^2\right)^{\frac{1}{2}
+\frac{1}{k}}\Psi_{M,F}^{1,2}\left(X^{
\frac{1}{k}}z^{\frac{2}{k}}u, X^{\frac{1}{2}}zv\right).
\end{align}
\end{proof}

\begin{lemma}\label{FT_formula}
The
function $\Psi_{M,F}$
satisfies, for all $i_1, i_2 \geq 0$,
\begin{align*}
&D_1^{i_1}D_2^{i_2} \Psi_{M,F}(x,y) \\&\ll_{i_1,i_2}
\left(M^{\frac{k-1}{2}} + |F|
M^{\frac{k}{4}-1}\right)^{i_1}\left(M^{\frac{k}{4}} +
\frac{|F|}{\sqrt{M}}\right)^{i_2} \|\phi\|_{C^{i_1+i_2}}
\|\psi\|_{C^{i_1 + i_2}},
\end{align*}
and therefore, for $u, v \neq 0$,
\begin{align}\label{FT_u_v_bound_k_finite}
&\Psi_{M,F}^{1,2}(u,v|\psi_j)\ll_{i_1,i_2}\\\notag&
M^{-\frac{k-2}{4}}\left(\frac{M^{\frac{k-1}{2}} + |F|
M^{\frac{k}{4}-1}}{|u|}\right)^{i_1}\left(\frac{M^{\frac{k}{4}
} + \frac{|F|}{\sqrt{M}}}{|v|}\right)^{i_2}\|\phi\|_{C^{
i_1+i_2}}\|\psi\|_{C^{i_1+i_2}}.
\end{align}
In terms of the frequency $F$, for $F \neq 0$,
\begin{equation}\label{FT_F_bound}
  \left|\Psi_{M,F}^{1,2}(u,0)\right| \ll \frac{M^{\frac{k}{4} +
\frac{1}{2}}}{|F|}
\|\Psi_{M,0}\|_1.
\end{equation}

\end{lemma}
\begin{proof}

The bounds on the derivatives are straightforward from the
observation
$x \ll \sqrt{M}$ and $y \ll x^{\frac{k}{2}}$, and the bound on the
Fourier
transform is deduced by integration by parts.

To prove the bound (\ref{FT_F_bound}), integrate
(\ref{Phi_FT}) by parts with respect to $y$ (note that $v =
0$) and use the bounds $x \ll \sqrt{M}$, $y \ll
M^{\frac{k}{4}}$.
\end{proof}

\begin{lemma}\label{mellin}
Let $\delta > 0$ and $\psi \in C^\infty(\bR^+)$ supported in $[\delta,
\infty]$,
satisfying, for all $a, j \geq 0$, $D^j (\psi(x) -1)x^a \to 0$ as $x \to
\infty$.  Set
\[ H\left(z\right) =
\Psi_{z^{-1},0}^{1,2}\left(0,0|\psi\right).\] For $s
\neq 0$ and $\frac{2}{k} s + \frac{1}{2} + \frac{1}{k} \neq -n$, $n \in
\zed_{\geq 0}$,
\[\tilde{H}\left(s\right) =
\frac{1}{k}\frac{\Gamma\left(\frac{1}{2}
\right)\Gamma\left(\frac{1}{2}-\frac{1}{k} +
  \frac{2}{k}s\right)}{\Gamma\left(1-\frac{1}{k} +
\frac{2}{k}s\right)}
\tilde{\phi}\left(\frac{1}{2} + \frac{1}{k} +
\frac{k-2}{k}s\right)
\tilde{\psi}\left(-s\right).\]
\end{lemma}

\begin{proof}
For $\Re s > 0$,
\begin{align*}
&\tilde{H}\left(s\right)
= 2
\int_0^\infty
\int_0^\infty \int_0^\infty \phi\left(x^k
-y^2\right)\psi\left(\frac{x^k-y^2}{x^2z}\right)
xyz^s \frac{dz}{z}\frac{dy}{y}\frac{dx}{x}
\\&= \int_0^\infty \int_0^1 \int_0^\infty
\phi\left(x^k\left(1-y\right)\right)\psi\left(x^{k-2}
\left(1-y\right)z\right)
x^{1 + \frac{k}{2}}y^{-\frac{1}{2}}z^{-s} \frac{dz}{z} dy
\frac{dx}{x}
\\&=\frac{1}{k}\int_0^\infty \phi\left(x\right)
x^{\frac{1}{2} + \frac{1}{k} + \frac{k-2}{k}s}\frac{dx}{x}
 \int_0^1 \left(1-y\right)^{\frac{2s}{k}
 - \frac{1}{2}-\frac{1}{k}}y^{-\frac{1}{2}} dy \int_0^\infty
\psi\left(z\right)z^{-s}\frac{dz}{z}
\\&=
\frac{1}{k}\frac{\Gamma\left(\frac{1}{2}
\right)\Gamma\left(\frac{1}{2}-\frac{1}{k} +
  \frac{2}{k}s\right)}{\Gamma\left(1-\frac{1}{k} +
\frac{2}{k}s\right)}
\tilde{\phi}\left(\frac{1}{2} + \frac{1}{k} +
\frac{k-2}{k}s\right)
\tilde{\psi}\left(-s\right).
\end{align*}
The conditions on $\psi$ guarantee that $\tilde{\psi}$
extends to a meromorphic
function, with a single simple pole of residue $-1$ at $s =
0$.  The formula
thus holds for $s$ not equal to a pole of the right hand
side, by analytic
continuation.
\end{proof}

\section{Proof of Theorems}
 The initial steps in the proofs of Theorems
\ref{quantitative_3_torsion_poincare} and
\ref{k_torsion_poincare} are made together, and
then the argument splits depending on $k = 3$ or $k > 3$, and
$f = 0$ or $f \neq
0$ when it is necessary to choose parameters.

The sums which appear in Theorems
\ref{quantitative_3_torsion_poincare} and
\ref{k_torsion_poincare} may be written
\begin{equation}\label{module_notation}
 \sum_{\substack{d \equiv 2 \bmod 4\\ \text{square-free}}}
\phi\left(\frac{d}{X}\right)
\sum_{\substack{\fa  \in P_k(-4d) \\ \fa = [a, b + \sqrt{-d}]}}
e\left(\frac{fb}{a}\right)\psi\left(\frac{\sqrt{d}}{Ta
}\right).
\end{equation}
To introduce the parameterization, let $H_k(-D)^*$ denote those classes in the
class group whose $k$th power is principal, and
write
\begin{equation}\label{finally_massaged}
 \mathscr{S}_{X,Y, f}= \sum_{\substack{d \equiv 2 \bmod 4\\ \text{square-free}}}
\phi\left(\frac{d}{X}\right) \sum_{\substack{(1) \neq  \fa
\text{  primitive},\\
[\fa]^k = [(1)] \in H(-4d)\\ \fa = [a, b + \sqrt{-d}]}}
e\left(\frac{fb}{a}
\right)\psi\left(\frac{\sqrt{d}}{Ta}\right).
\end{equation}
This sum counts ideals from classes that are order dividing $k$, but those of
order less than $k$ do not appear due to the conditions which are imposed
upon $T$ and the support of $\psi$.

 In the case that $\psi \equiv 1$ near $\infty$ it is
convenient to
localize further, so as to consider Heegner points having
imaginary part in
dyadic intervals.  Let $\sigma^\times$ be the smooth
 multiplicative partition of unity function of the previous
section (see (\ref{sigma_properties})), so that  $\sum_{j \in
\zed}
\sigma^\times(e^jx) = 1$.   Define, for $j \in \zed$, and
real $x$
\[
  Y_j = e^j, \qquad  \psi_j(x) = \left\{\begin{array}{lll}\psi\left(\frac{Y_j
x^{\frac{1}{2}}}{T}\right)
\sigma^\times\left(x^{\frac{1}{2}}\right)&& x > 0\\ 0 && x \leq 0
\end{array}\right..
\]
If $j < \log T - 12$ the support condition gives $\psi_j \equiv 0$.
We may now write
\begin{align}
 \notag  \mathscr{S}_{X,Y, f}&= \sum_{j } \sum_{\substack{d
\equiv 2 \bmod 4\\ \text{square-free}}}
\phi\left(\frac{d}{X}\right) \sum_{\substack{(1) \neq  \fa
\text{ primitive},\\
[\fa]^k = [(1)] \in H(-4d)\\ \fa = [a, b + \sqrt{-d}]}}
e\left(\frac{fb}{a}\right)\psi_j\left(
\frac{d}{Y_j^2 a^2}\right)
 \\ \label{sum_decomp} & =  \sum_{\log T - 12 \leq j \ll \log X }
\mathscr{S}_{j}.
\end{align}

\subsection{Global parameterization}
Recall the definition from (\ref{def_Phi_X_Y_a_f}),
\[
 \Phi_{X, Y,f}\left(x,y,z| \psi\right)
=
\phi\left(\frac{x^k -
y^2}{Xz^2}\right)\psi\left(\frac{x^k
  -y^2}{Y^2 x^2z^2}\right)e\left(-\frac{f y}{xz}\right).
\]
 Solving $d = \frac{\ell m^k - \ell^2n^2}{t^2},$ and observing
the property of
fractions
\[
  \frac{n t^{-1}}{m} \equiv \frac{-n m^{-1}}{t} +
\frac{n}{mt} \bmod 1,\] where $
tt^{-1}\equiv 1 \bmod m,$ and $mm^{-1}\equiv 1  \bmod t
$,
write
\begin{align}\label{def_S_j}
 \mathscr{S}_{j } &= \sum_{\substack{\ell,m,t \in \zed^+, n
\in \zed\\ \cC_1}} e\left(\frac{fnm^{-1}}{t}\right)\Phi_{X,Y_j,
f}\left(\ell m,
\ell^{\frac{k+1}{2}}n, \ell^{\frac{k-1}{2}}t\Big| \psi_j \right)\\& \notag
\qquad \times
\sum_{s^2|\frac{\ell m^k-\ell^2n^2}{t^2}} \mu(s).
\end{align}
In this expression, $\cC_1$ indicates the local conditions
\[
\cC_1 = \left\{\ell\; \square\text{-free}, (\ell mn,t)=(2\ell, m)=1, \ell m^k
-\ell^2n^2 \equiv 2t^2 \bmod 4t^2 \right..
\]

Note that support of $\psi_j$ and $\phi$ imposes the following restrictions on
the summations variables
\begin{align} \label{length_lt_sum} \ell m &\ll \frac{\sqrt{X}}{Y}\\ \notag
\ell^{\frac{k-1}{2}}t &\ll X^{\frac{k-2}{4}
}Y^{-\frac{k}{2}}.\end{align}

 Splitting the sum over $s$ at
parameter $Z$
write $\mathscr{S}_j = \mathscr{M}_{j} +  \mathscr{E}_j$
as a main term plus an error term.  In the main term, perform M\"{o}bius
inversion with variable $s_1$ to eliminate the
co-primality condition between $\ell$ and $m$. Write $s_1 \ell := \ell$, $s_1
m:= m$.  Thus\footnote{The subscripts $X,
Y_j, f$ are suppressed.}
\begin{align}\label{def_main_term}
 &\cM_{j} =  \sum_{s < Z, s_1} \mu(s s_1) \\\notag
&  \times\sum_{\substack{\ell,m,t \in \zed^+, n
\in \zed\\ \cC_2 }}
e\left(\frac{fns_1^{-1}m^{-1}}{t}\right)\Phi\left(s_1^2\ell m,
(s_1 \ell)^{\frac{k+1}{2}}n,
(s_1 \ell)^{\frac{k-1}{2}}t\Big| \psi_j\right);\\
\notag &\cC_2  = \left\{\begin{array}{l}\ell\;
\text{square-free},\\ (s_1\ell mn,t) =(\ell,ss_1)= 1,  \\
s_1^{k+1}\ell m^k - s_1^2\ell^2n^2
\equiv 2 s^2t^2 \bmod 4s^2t^2 \end{array} \right.
\end{align}
and
\begin{align*}
 \mathscr{E}_{j} &= \sum_{\substack{\ell,m,t \in \zed^+, n
\in \zed\\ \cC_3}}
e\left(\frac{fnm^{-1}}{t}\right)\Phi
\left(\ell m,
\ell^{\frac{k+1}{2}}n, \ell^{\frac{k-1}{2}}t\Big| \psi_j \right)
\sum_{\substack{s^2|\frac{\ell
m^k-\ell^2n^2}{t^2}\\ s \geq  Z}} \mu(s);\\
\cC_3 &= \left\{\begin{array}{l} \ell\; \text{square-free}\\ (\ell mn,t) =
(\ell, m)
= 1\\  \ell m^k  - \ell^2n^2
\equiv 2 t^2 \bmod
4t^2\end{array} \right..
\end{align*}
The next section is concerned with proving the following
evaluation of the main term.

\begin{proposition}\label{main_term_prop}
Let $k \geq 3$, odd, and let $c_k$ be the constant of Theorem
\ref{k_torsion_secondary_term}.  In the case $f = 0$, in the range $Z \ll
T^{\frac{k}{4}}X^{\frac{1}{2} - \frac{k}{8} -\epsilon}$,
$\cM = \sum \cM_j$ satisfies
\begin{align*}
 \cM =& \frac{6}{\pi^3} \tilde{\phi}(1)\tilde{\psi}(-1)
\frac{X}{T}+
\psi(\infty)\tilde{\phi}\left(\frac{1}{2} + \frac{1}{k}\right) c_k
X^{\frac{1}{2}
+\frac{1}{k}}\\ &
+ O\left(X^{\frac{1}{2} +
\frac{1}{(2k-2)} + \epsilon}\right) +
O\left(X^{1+\epsilon}T^{-1}Z^{-1}\right) +
O\left(X^{\frac{k}{4} + \epsilon}T^{-\frac{k}{2}}\right).
\end{align*}

When $f \neq 0$, for  $Z \ll |f|^{-\frac{1}{2}}
X^{\frac{-k}{8}
+ \frac{1}{2}
-\epsilon}T^{\frac{k}{4}-\frac{1}{2}},
$
\[
 \mathscr{M}_f = O\left(X^{\frac{k}{4}
+
\epsilon}T^{-\frac{k}{2}}\right) + O\left(X^{\frac{1}{2} + \frac{1}{k}
+ \epsilon}\right) +
\delta_{k=3} O\left(X^{\frac{7}{8} +
\epsilon}T^{-\frac{1}{4}}\right).
\]
\end{proposition}
In the final section, Section \ref{sieve_section}, the sieving error
term is estimated.
\begin{proposition}\label{error_term_prop}
We have
\[
 \mathscr{E} = \sum_j \mathscr{E}_j \ll \frac{X^{1 +
\epsilon}}{TZ} +
\frac{X^{\frac{k}{4} + \epsilon}}{T^{\frac{k}{2}}}.
\]
\end{proposition}
One easily obtains by Mellin inversion
\begin{equation}\label{square_free_count}
 \sum_{\substack{d \equiv 2 \bmod 4\\ d
\;\square\text{-free}}}
\phi\left(\frac{d}{X}\right) =
\frac{2}{\pi^2}\tilde{\phi}(1)X +
O\left(X^{\frac{1}{2}}\right).
\end{equation}
The deductions Theorems
\ref{quantitative_3_torsion_poincare} and
\ref{k_torsion_poincare} are as follows.

\begin{proof}[Proof of Theorem
\ref{quantitative_3_torsion_poincare}]
Recall that this treats the case $k = 3$.

When $f = 0$ choose $Z = T^{\frac{3}{4}} X^{\frac{1}{8}
-\epsilon}$
to obtain the asymptotic of the Theorem
 with error term bounded by $O\left(\frac{X^{\frac{7}{8} +
\epsilon}}{T^{\frac{7}{4}}}\right) + O\left(X^{\frac{3}{4}
+\epsilon}\right).$

When $f \neq 0$ choose $Z = |f|^{-\frac{1}{2}}
T^{\frac{1}{4}}X^{\frac{1}{8}  -\epsilon}$ to obtain
the bound
$
 O\left(|f|^{\frac{1}{2}}
\frac{X^{\frac{7}{8}}}{T^{\frac{5}{4}}}\right) +
O\left(X^{\frac{5}{6} + \epsilon}\right)
$
as required.
\end{proof}

\begin{proof}[Proof of Theorem \ref{k_torsion_poincare}]
 When $f = 0$, choose $Z = T^{\frac{k}{4}}X^{\frac{1}{2} -
\frac{k}{8}-\epsilon}$ to
obtain the asymptotic of the Theorem with error terms of
size
$
 O\left(\frac{X^{\frac{k}{4}+
\epsilon}}{T^{\frac{k}{2}}}\right) +
O\left(X^{\frac{1}{2} + \frac{1}{2k-2}+\epsilon}\right).
$

When $f \neq 0$, choose $Z =  |f|^{-\frac{1}{2}}
T^{\frac{k}{4}-\frac{1}{2}}X^{\frac{-k}{8} + \frac{1}{2}
-\epsilon}$ to obtain
a bound of
$
 O\left(\frac{X^{\frac{k}{4}+\epsilon}}{T^{\frac{k}{2}}}
\right) +
O\left(\frac{|f|^{\frac{1}{2}}X^{\frac{k}{8} +
\frac{1}{2}+\epsilon}}{T^{\frac{k}{4} + \frac{1}{2}}}\right)
+
O\left(X^{\frac{1}{2} + \frac{1}{k} + \epsilon}\right).
$
\end{proof}

\subsection{Evaluation of main term}
Control the local conditions in $\mathscr{M}_{j}$ by setting $s_2
= (s, t)$ and $s_3 = \GCD(s, m, n)$.  Then replace $s_2 t:= t$, $s_3m:= m$,
$s_3 n:= n$.  Thus\footnote{In this section the indices $X, Y, f$
are suppressed}
\begin{align*}
\cM_j =&\sum_{\substack{s < Z, s_1\\ s = s_2s_3s_4\\ \text{odd}}}\mu(ss_1)
\mathop{\sum\sum}_{\substack{(\ell, t) \in (\zed^+)^2, (m,n) \in
\zed^2\\ \cC_4}}
\\& e\left(\frac{fn
s_1^{-1}m^{-1}}{s_2 t}\right)
\Phi\left(s_1^2 s_3 \ell
m, (s_1\ell)^{\frac{k+1}{2}} s_3 n, (s_1 \ell)^{\frac{k-1}{2}}s_2 t \right);\\
\cC_4 =& \left\{\begin{array}{l} \ell \;
    \text{square-free}\\  (t, s_1s_3s_4) =(\ell, s s_1t) =
1\\m\text{ odd} \\ (mn, s_2s_4t)=1
 \\
    \ell s_1^{k+1}s_3^km^k - \ell^2 s_1^2 s_3^2 n^2
\equiv
    2s_2^4s_4^2t^2 \bmod{4s_2^4s_4^2t^2} \end{array}\right..
\end{align*}
Set $\cM_j = \cM_{j,e} + \cM_{j,o}$ according as $\ell$ is even or
odd. When $\ell$ is even the
condition
at 2 is guaranteed so that on replacing $\ell$ by $\frac{\ell}{2}$,
\begin{align*} \cM_{j,e}&= \sum_{\substack{s
<
Z, s_1  \\s = s_2s_3s_4\\ \text{odd}}}\mu(ss_1)
\mathop{\sum\sum}_{\substack{(\ell, t) \in (\zed^+)^2, (m,n) \in
\zed^2\\\cC_{5,e}}}\\& \qquad e\left(\frac{f
s_1^{-1}n m^{-1}}{s_2 t}\right)
\Phi\left(2s_1^2s_3 \ell
m,(2 s_1 \ell)^{\frac{k+1}{2}} s_3  n, (2 s_1\ell)^{\frac{k-1}{2}}s_2 t\right)\\
\cC_{5,e}&= \left\{ \begin{array}{l} \ell \;
    \text{square-free}\\ (\ell, 2s_1st) = (t, 2s_1s_3s_4) =
1 \\
m\text{  odd}\\ (mn, s_2s_4t)=1
\\
    (2s_1^2s_3 \ell m)^k\equiv
\left((2s_1\ell)^{\frac{k+1}{2}}s_3 n\right)^2
\bmod{s_2^4s_4^2t^2} \end{array}\right..
\end{align*}
Setting apart the sum over $m$ and $n$, write
\begin{align*}
\cM_{j,e}&= \sum_{\substack{s = s_2s_3s_4 <
Z, s_1 \\ \text{odd}}}\mu(ss_1)
\sum_{\substack{(\ell, t) \in (\zed^+)^2\\ \ell \;
    \square\text{-free}\\ (\ell, 2s_1st) = (t, 2s_1s_3s_4) =
1}} \cM_{j,e, s,\ell,
t}.
\end{align*}
When $\ell$ is odd,
\begin{align*}\cM_{j,o} &= \sum_{\substack{s
< Z, s_1 \\s = s_1s_2s_3s_4\\
\text{odd}}}\mu(ss_1)
\mathop{\sum\sum}_{\substack{(\ell, t) \in (\zed^+)^2, (m,n) \in \zed^2 \\
\cC_{5,o} }}
\\&\qquad  e\left(\frac{f
s_1^{-1}nm^{-1}}{s_2 t}\right)\Phi\left(s_1^2s_3 \ell m, (s_1
\ell)^{\frac{k+1}{2}}
s_3  n,
(s_1\ell)^{\frac{k-1}{2}} s_2t \right)
\\
\cC_{5,o} &= \left\{ \begin{array}{l} \ell
\;\text{square-free}\\  (t, s_1s_3s_4) =(\ell, 2s_1st) =
1\\  m \text{ odd} \\ (mn, s_2s_4t)=1 \\
     (s_1^2s_3 \ell m)^k -\left((s_1 \ell)^{\frac{k+1}{2}}
s_3 n\right)^2
\equiv 2s_2^4s_4^2t^2
\bmod{4s_2^4s_4^2t^2} \end{array}\right..
\end{align*}
As above, write
\begin{align*}
\cM_{j,o}&= \sum_{\substack{s = s_2s_3s_4 < Z, s_1 \\
\text{odd}}}\mu(ss_1)
\sum_{\substack{(\ell, t) \in (\zed^+)^2\\ \ell \;
    \square\text{-free, odd}\\ (\ell, s_1st) = (t, s_1s_3s_4) =
1}}\cM_{j,o,
s,\ell,t}.
\end{align*}
We show the analysis in the even case.  The odd case may be
handled similarly.

\subsection{Local parameterization}
By Proposition \ref{local_parametrization} the sum
over $(m,n)$ in $\cM_{j,e}$ is parametrized by setting
\begin{align*}
 2s_1^2s_3 \ell m &= (2 s_1s_3 \ell w)^2 + (2a + 1) \cdot 2
 s_1^2s_3\ell\cdot N\\
  (2 s_1\ell)^{\frac{k+1}{2}}s_3 n &= (2  s_1s_3\ell w)^k + (2a
+ 1) \cdot
\frac{k}{2}  (2 s_1 s_3 \ell w)^{k-2} (2s_1^2 s_3 \ell)\cdot
N\\ &
\qquad  + b\cdot (2 s_1\ell)^{\frac{k+1}{2}}s_3 \cdot
N^2,
\end{align*}
where
\begin{equation}\label{ranges}
a, b \in \zed, \qquad w \in
\left(\zed /N \zed\right)^\times, \qquad N = s_2^2 s_4 t.
\end{equation}
Thus\footnote{In this section we abbreviate $\cM =\cM_{ j,e, s,
\ell, t}$.}
\begin{align*} \cM =&
\frac{1}{s_2^2 s_4 t}\sum_{\substack{0 \leq w <s_2^4 s_4^2 t^2\\ (w,
s_2s_4t)=1}} \sum_{\substack{a, b \in
\zed}}e\left(\frac{\tilde{f}w^{k-2}}{s_2t}\right)
 \Phi(A + Ba, C + Da + Eb, z)\end{align*}
where
\begin{align*}
 A &= (2\ell s_1s_3 w)^2 + 2 \ell
 s_1^2s_3\cdot s_2^2 s_4 t\\
 B &= 4 \ell s_1^2s_3\cdot s_2^2 s_4 t\\
 C &= (2\ell s_1s_3 w)^k + k\cdot s_1^2 s_3 \ell(2s_1s_3\ell w)^{k-2}\cdot
s_2^2 s_4 t\\
 D &= 2k\cdot s_1^2 s_3  \ell(2s_1s_3\ell w)^{k-2} \cdot s_2^2s_4 t\\
 E &= (2s_1 \ell)^{\frac{k+1}{2}}s_3\cdot s_2^4s_4^2t^2\\
 \tilde{f} &= f \cdot
2^{\frac{k-3}{2}}s_1^{\frac{k-3}{2}}s_3^{k-2}\ell^{\frac{k-3}{2}} \\
 z &= (2s_1\ell)^{\frac{k-1}{2}}s_2t.
\end{align*}

\begin{lemma}Keep the definitions of $A-E, \tilde{f}, z$ above, and set
 \begin{align}\label{definition_of_M_F} M &=
\frac{X^{1-\frac{2}{k}}}{Y_j^2(2s_1\ell)^{2-\frac{2}{k}}
(s_2t)^{\frac{4}{k}}}, \qquad F= \frac{f
X^{\frac{1}{2}-\frac{1}{k}}}{(2s_1\ell)^{1-\frac{1}{k}}(s_2
t)^{\frac{2}{k}}}.\end{align}
Define
\begin{align*}
 U_0 & = \frac{s X^{\frac{k-3}{2}}}{Y_j^{k-1} (s_1
\ell)^{k-3}s_2 t }, \qquad
 \forall f \neq 0, \; U_f = \frac{|f| s X^{\frac{k}{4}-1
}}{Y_j^{\frac{k}{2} -2}(s_1\ell)^{\frac{k}{4}-1}}.
\end{align*}
Subject to the constraint on $Z$
\begin{equation*}
 Z \leq \left\{ \begin{array}{lll}
 X^{\frac{1}{2} - \frac{k}{8}
-\epsilon}T^{\frac{k}{4}} && f = 0 \\
|f|^{-\frac{1}{2}} X^{\frac{1}{2} - \frac{k}{8}
-\epsilon}T^{\frac{k}{4}-\frac{1}{2}}&& f \neq 0
 \end{array}\right.,
\end{equation*}
for any $N>0$,
\begin{align*}
  \cM&=  O_{N}\left(X^{-N}\right) + \Delta_f + E_f \\
  \Delta_f &= \frac{(z^2X)^{\frac{1}{2} +
\frac{1}{k}}}{BE}  S_{k-2}(0, \tilde{f}s_2s_4t;
s_2^2s_4t) \Psi_{M,F}(0,0)\\
E_f &=
\frac{(z^2X)^{\frac{1}{2} +
\frac{1}{k}}}{BE}\\& \qquad \times\sum_{0 \neq |u| \leq U_f
X^{\epsilon}}(-1)^u S_{k-2}(\ell s_3 u,
\tilde{f}
s_2 s_4; s_2^2 s_4 t)\Psi^{1,2}_{M, F}\Biggl(\frac{\left(z^2 X
\right)^{\frac{1}{k}}u}{B}, 0\Biggr).
 \end{align*}

\end{lemma}

\begin{proof}
Applying Poisson summation in the $a$ and $b$ variables, and
evaluating the Fourier
transform by
applying Lemma \ref{fourier_shift},
 \begin{align*}&\cM =
\frac{1}{BE} \frac{1}{s_2^2 s_4 t}\sum_{\substack{0 \leq w <s_2^4 s_4^2 t^2\\
(w,
s_2s_4t)=1}}\sum_{\substack{u,v \in
\zed}}e\left(\frac{\tilde{f}w^{k-2}}{
s_2t}+\frac{Au}{B} + \left(\frac{BC - AD}{BE} \right) v\right)\\&
\quad\times
  \Phi^{1,2}\left(\frac{u}{B} - \frac{Dv}{BE},
\frac{v}{E}, z\right).
\end{align*}

Applying Lemma \ref{FT_formula},
\begin{align}\label{v_sum}
 &\cM=\frac{(z^2X)^{\frac{1}{2} + \frac{1}{k}}}{BEs_2^2s_4t}
\sum_{\substack{0 \leq w < s_2^4 s_4^2 t^2\\ (w, s_2s_4t)=1}}
\sum_{\substack{u,v \in
    \zed
}}e\left(\frac{\tilde{f}w^{k-2}}{s_2 t}+\frac{Au}{B} + \left(\frac{BC - AD}{BE}
\right) v\right)
 \\& \notag \times \Psi^{1,2}_{M, F}\Biggl(\frac{\left(z^2 X
\right)^{\frac{1}{k}}}{B}
\left(u - \frac{Dv}{E} \right),
  \frac{zX^{\frac{1}{2}}v}{E}\Biggr).
\end{align}

Decay of the Fourier transform is now used to truncate the ranges
of summation. By rapid decay of $\Psi_{M,F}^{1,2}$ in the first and second
slots
((\ref{FT_u_v_bound_k_finite}) of Lemma \ref{FT_formula}), the sums
over $u, v$ and $w$
above are bounded in length by polynomials in $X$, with
negligible error.

We first argue that we may discard all terms with $v \neq 0$ with negligible
error.  By decay in the second slot, those terms satisfying
\[
\left[ M^{\frac{k}{4}} + \frac{|F|}{\sqrt{M}}\right]\left(
\frac{s_1s_2^3s_3s_4^2
\ell t}{X^{\frac{1}{2}}}\right) < X^{-\epsilon},  \epsilon > 0
\]
are bounded by, for all $N>0$, $O_{N}\left(X^{-N}\right)$.
Suppose first that $f = 0$ so that $F = 0$.  Then, using (\ref{length_lt_sum})
\begin{align*}
s_1^{\frac{k+1}{2}} s_2 \ell^{\frac{k-1}{2}} t &\ll
X^{\frac{k-2}{4}}Y^{-\frac{k}{2}}, \qquad
M =
\frac{X^{1-\frac{2}{k}}}{Y_j^2(2s_1\ell)^{2-\frac{2}{k}}
(s_2t)^{\frac{4}{k}}}
\end{align*}
we have
\[
M^{\frac{k}{4}} \frac{s_1s_2^3 s_3 s_4^2 \ell t}{X^{\frac{1}{2}}}
\ll
\frac{X^{\frac{k}{4}-1}s_2^2s_3s_4^2}{Y_j^\frac{k}{2}} \ll
X^{\frac{k}{4} -1}
 T^{-\frac{k}{2}} Z^2
\]
and so the condition $
 Z \ll  X^{\frac{1}{2} - \frac{k}{8}
-\epsilon}T^{\frac{k}{4}}$ suffices.

When $f \neq 0$, one must consider in addition,
\[
  \frac{|F|}{\sqrt{M}} \frac{s_1s_2^3 s_3 s_4^2 \ell
t}{X^{\frac{1}{2}}} = |f|
s_1s_2^3 s_3 s_4^2 \ell t X^{\frac{-1}{2}}Y_j \ll |f|s_2^2 s_3 s_4^2
X^{\frac{k}{4}-1} Y_j^{-\frac{k}{2}+1},
\]
where in the last inequality we again use (\ref{length_lt_sum}).
Therefore, for $f \neq 0$ the condition $
  Z \ll |f|^{-\frac{1}{2}} X^{\frac{1}{2} - \frac{k}{8}
-\epsilon}T^{\frac{k}{4}-\frac{1}{2}}$ suffices.

Thus in the given ranges for $Z$ we may assume that $v = 0$ and now truncate
the
sum
over $u$.  This is
negligible beyond the range
 \[
  |u| \ll X^\epsilon \left[ M^{\frac{k-1}{2}} +
|F|M^{\frac{k}{4}-1}\right]
\frac{s_1^{1 + \frac{1}{k}} s_2^{2-\frac{2}{k}}s_3s_4
\ell^{\frac{1}{k}}t^{1-\frac{2}{k}}}{X^{\frac{1}{k}}}.
\]
When $f = 0$, this gives the restriction
\begin{equation}\label{f_0_u_length}
  |u| \leq   U_0 X^{\epsilon}.
\end{equation}
 When $f \neq 0$ the second term dominates, and we have the
restriction
\begin{equation}\label{f_not_0_u_length}
 |u| \leq
U_fX^\epsilon.
\end{equation}

With $v = 0$, the inner sum over $w$ in (\ref{v_sum})
becomes the complete sum
\begin{align*}
(-1)^u s_2^2 s_4 t \cdot S_{k-2}\left(  \ell s_3 u, \tilde{f}s_2
s_4; s_2^2 s_4 t\right),
\end{align*}
completing the evaluation.
\end{proof}

\subsubsection{Evaluation of the diagonal
$\Delta_0$}\label{diagonal_section}
 When $f = 0$, $\Delta_0$ is a diagonal main term contribution.
Write
$\Delta_{0,j,e}$ to indicated $\Delta_0$ for the even terms attached to
$\cM_j$.  Since $S_{k-2}(0,0; s_2^2 s_4 t) =
\varphi(s_2^2 s_4 t)$ we have
\[
 \Delta_{0,j,e} = \sum_{\substack{s = s_2s_3s_4< Z, s_1 \\
\text{odd}}}\mu(ss_1)
\sum_{\substack{(\ell, t) \in (\zed^+)^2\\ \ell \;
    \square\text{-free}\\ (\ell, 2st) = (t, 2s_1s_3s_4) =
1}} \frac{X^{\frac{1}{2} + \frac{1}{k}}\varphi(s_2^2s_4t)
\Psi_{M_j,0}^{1,2}(0,0)}{2^{2+\frac{1}{k}}s_1^{2 +
\frac{1}{k}}s_2^{5-\frac{2}{k}}s_3^2s_4^3 \ell^{1 +
\frac{1}{k}}t^{2-\frac{2}{k}}}
\]
 where $M_j = \frac{X^{1-\frac{2}{k}}}{Y_j^2(2s_1\ell)^{2 -
\frac{2}{k}}(s_2
t)^{\frac{4}{k}}}.$  As a first step we remove the
restriction $s > Z$.
It follows from (\ref{L_1_bound}) that
$\Psi_{M_j, 0}^{1,1}(0,0)
\ll X^\epsilon M_j^{-\frac{k-2}{4}}$.
Substituting this bound,  the sum over $s \geq
Z$ is bounded by (use $Y \ll X^{\frac{1}{2} - \frac{1}{k}}$ in bounding the sum
over $s_1\ell$)
\[
  \sum_{s=s_2s_3s_4 \geq
Z}\sum_{\substack{(s_1\ell)^{\frac{k-1}{2}} s_2t \\ \ll
X^{\frac{k-2}{4}} Y_j^{\frac{-k}{2}}}} \frac{X^{\frac{3}{2} - \frac{k}{4}
+
\epsilon}Y_j^{\frac{k-2}{2}}}{
  s_1^{\frac{7}{2}-\frac{k}{2}} s_2^2 s_3^2 s_4^2
\ell^{\frac{5}{2}-\frac{k}{2}}} \ll X^{1+
\epsilon}Y_j^{-1}Z^{-1}.
\]

 Next remove the partition of unity.  Recall that
$Y_j = e^j$, and that $\psi_j$ is supported on $x>0$, defined there by
\[
 \psi_j(x) = \psi\left(\frac{Y_j x^{\frac{1}{2}}}{T}\right)
\sigma^\times\left(x^{\frac{1}{2}}\right).
\]
Put $\psi_{\frac{1}{2}}(x) = \psi\left(\frac{x^{\frac{1}{2}}}{T}\right)$ for $x
> 0$, $\psi_{\frac{1}{2}}(x) = 0$ for $x \leq 0$.

\begin{lemma}
  For arbitrary $M > 0$ we have the equality
 \[\sum_{j}
 \Psi_{\frac{M}{Y_j^2},0}^{1,2}\left(0,0|\psi_j\right)=
\Psi_{M,0}^{1,2}\left(0,0|\psi_{\frac{1}{2}}\right).\]
\end{lemma}

\begin{proof}
The left hand side is
 \begin{align*}&\int_{\bR^2}\phi\left(x^k - y^2\right)\left[
\sum_{j }
 \psi_j\left(\left(x^k-y^2\right) \frac{M}{Y_j^2 x^2} \right)\right]
dx dy= \Psi_{M,0}^{1,2}\left(0,0|\psi_{\frac{1}{2}}\right). \end{align*}
\end{proof}
Applying the lemma,
\begin{align*}
\Delta_{0,e} &= \sum_{j}
\Delta_{0,j,e} =O(X^{1+\epsilon} T^{-1}Z^{-1})\\&+
 X^{\frac{1}{2} + \frac{1}{k}} \sum_{s = s_1
s_2s_3s_4}\mu(s)\sum_{\substack{(\ell, t) \in (\zed^+)^2\\ \ell \;
     \square\text{-free}\\ (\ell, 2st) = (t, 2s_1s_3s_4) =
1}}
 \frac{\varphi(s_2^2s_4t)
\Psi_{\ast,0}^{1,2}\left(0,0\big|\psi_{\frac{1}{2}}\right)}{2^{2 +
\frac{1}{k}}s_1^{2 + \frac{1}{k}}s_2^{5-\frac{2}{k}}s_3^2s_4^3
\ell^{1+\frac{1}{k}}t^{2-\frac{2}{k}}}
\end{align*}
 where $\ast$ stands in for
$\frac{X^{1-\frac{2}{k}}}{(2s_1\ell)^{2 -
\frac{2}{k}}(s_2 t)^{\frac{4}{k}}}$.

  Dropping the error, we now evaluate the main term by
Mellin inversion, using
the formula of Lemma \ref{mellin} for the Mellin transform
of
$\Psi_{z^{-1},0}^{1,2}(0,0)$.  This yields the main term as the integral
\begin{align*}&\frac{\Gamma\left(\frac{1}{2}\right)}{k}\frac{X^{\frac{1}{2} +
\frac{1}{k}}}{2^{2+
\frac{1}{k}}} \\&
\oint_{(2)}\frac{\Gamma\left(\frac{1}{2}- \frac{1}{k} +
\frac{2\alpha}{k}
\right)}{\Gamma\left(1 - \frac{1}{k} +
\frac{2\alpha}{k}\right)}\tilde{\phi}\left(\frac{1}{2} +
\frac{1}{k} +
\frac{k-2}{k} \alpha\right)
\tilde{\psi}_{\frac{1}{2}}\left(-\alpha\right)
\frac{X^{(1-\frac{2}{k})\alpha}}{2^{(2 -
\frac{2}{k})\alpha}} F(\alpha)d\alpha
\end{align*}
where $\tilde{\psi}_{\frac{1}{2}}(\alpha) =
2T^{2\alpha}\tilde{\psi}(2\alpha)$ and
\begin{align*}
 & F(\alpha)\\& = \sum_{\substack{s_1s_2s_3s_4 = s\\
\text{odd}}}
\frac{\mu(s)}{s_1^{2 + \frac{1}{k} + (2 -
\frac{2}{k})\alpha}s_2^{5 -
\frac{2}{k} + \frac{4\alpha}{k}} s_3^2
s_4^3}\sum_{\substack{(\ell, t) \in
(\zed^+)^2\\ \ell \;
     \square\text{-free}\\ (\ell, 2st)= \\ (t, 2s_1s_3s_4) =
1}}
\frac{\varphi(s_2^2 s_4 t)}{\ell^{1 + \frac{1}{k} + (2 -
\frac{2}{k})\alpha}t^{2
- \frac{2}{k} + \frac{4\alpha}{k}}}\\
&= \zeta\left(1 - \frac{2}{k} +
\frac{4\alpha}{k}\right) G(\alpha)
\end{align*}
 where
\begin{align*}
&  G(\alpha) = \left(1 - \frac{1}{2^{1 -\frac{2}{k} +
\frac{4\alpha}{k}}}\right)\\&\prod_{p \text{
odd}}\biggl[1 -
 \frac{2}{p^2} + \frac{1}{p^3} + \frac{1}{p^{1 + \frac{1}{k}
+
(2-\frac{2}{k})\alpha}}
  - \frac{1}{p^{2- \frac{1}{k} + (2 + \frac{2}{k})\alpha}}-
\frac{1}{p^{2 -
\frac{2}{k} + \frac{4\alpha}{k}}}\\&\qquad\qquad\qquad\qquad
+ \frac{1}{p^{3 -
\frac{2}{k} + \frac{4\alpha}{k}}}- \frac{1}{p^{2 +
\frac{1}{k} +
(2-\frac{2}{k})\alpha}} + \frac{1}{p^{3 - \frac{1}{k} + (2 +
\frac{2}{k})\alpha}}  \biggr].
\end{align*}
 $G$ is holomorphic in $\Re(\alpha)> \frac{-1}{2k-2}$.
Shifting the contour to $\Re(\alpha) = \frac{-1}{2k-2}
+ \epsilon$, we
pass a pole at $\frac{1}{2}$, and, depending on $\psi$,
possibly a second pole
at $\alpha = 0$.  We have
$G\left(\frac{1}{2}\right) = \frac{32}{\pi^4}$
and
\[
 G(0) = \frac{8}{\pi^2}\left(1 -
\frac{1}{2^{1-\frac{2}{k}}}\right)\prod_{p
\text{ odd}}\left[1 +
\frac{1}{p+1}\left(\frac{1}{p^{\frac{1}{k}}} -
\frac{1}{p^{1-\frac{2}{k}}} -\frac{1}{p^{1-\frac{1}{k}}}
-\frac{1}{p}\right)\right]
\]
Thus
\begin{align*}
  \Delta_{0,e} = &O\left(X^{\frac{1}{2}+\frac{1}{(2k-2)} + \epsilon}\right) +
O\left(X^{1+\epsilon}T^{-1}Z^{-1}\right)+\frac{2}{\pi^3}
\tilde{\phi}(1)\tilde{\psi}(-1) \frac{X}{T} \\&+
\psi(\infty)\tilde{\phi}\left(\frac{1}{2} + \frac{1}{k}\right)
\frac{1}{k\pi^{\frac{3}{2}}}
\frac{\Gamma(\frac{1}{2}-\frac{1}{k})}{\Gamma(1-\frac{1}{k})
}X^{\frac{1}{2}
+\frac{1}{k}}
\left(2^{1-\frac{1}{k}} -
2^{\frac{1}{k}}\right)\\ &\qquad\qquad \times \prod_{p \text{ odd}}\left[1 +
\frac{1}{p+1}\left(\frac{1}{p^{\frac{1}{k}}} -
\frac{1}{p^{1-\frac{2}{k}}}
-\frac{1}{p^{1-\frac{1}{k}}} -\frac{1}{p}\right)\right].
\end{align*}
The analysis of $\Delta_{0,o}$ is entirely analogous.  It
yields,
\begin{align*}
 \Delta_{0,o} = &O\left(X^{\frac{1}{2}+\frac{1}{(2k-2)} + \epsilon}\right) +
O\left(X^{1+\epsilon}T^{-1}Z^{-1}\right)+ \frac{4}{\pi^3}
\tilde{\phi}(1)\tilde{\psi}(-1) \frac{X}{T}\\& +
\psi(\infty)\tilde{\phi}\left(\frac{1}{2} + \frac{1}{k}\right)
\frac{1}{k\pi^{\frac{3}{2}}}
\frac{\Gamma(\frac{1}{2}-\frac{1}{k})}{\Gamma(1-\frac{1}{k})
}X^{\frac{1}{2}
+\frac{1}{k}} \\&\qquad\qquad\qquad\times  \prod_{p \text{
odd}}\left[1 +
\frac{1}{p+1}\left(\frac{1}{p^{\frac{1}{k}}} -
\frac{1}{p^{1-\frac{2}{k}}}
-\frac{1}{p^{1-\frac{1}{k}}} -\frac{1}{p}\right)\right].
\end{align*}
 Combining these two expressions together obtains the main
term of Proposition
\ref{main_term_prop}.

\subsubsection{Bound for the off-diagonal $f=0$, $u \neq
0$}\label{off_diagonal_section}
Write $E_{0,j,e,s, \ell, t}$ for the even terms associated to $E_{0}$ coming
from $\cM_j$.   It follows from Lemma \ref{S_f_u_lemma} that
\[\left|S_{k-2}(\ell s_3
u, 0; s_2^2 s_4 t)\right|\ll (u, s_2^2 s_4
t)^{\frac{1}{2}}(s_2^2s_4 t)^{\frac{1}{2} +
\epsilon}.\] Actually we could quite easily extract the sign
and get much more
cancellation, but anyway, this is not the limiting error
term.

In view of the restriction $u \ll U_0 X^\epsilon$ (see
(\ref{f_0_u_length})) we obtain
\begin{align*}
 E_{0,j,e, s, \ell, t} &\ll \frac{X^{\frac{1}{2} +
\frac{1}{k}}\|\Psi_{M,0}(\psi_j)\|_1}{s_1^{2+
\frac{1}{k}}s_2^{4 -
\frac{2}{k}}s_3^2s_4^{\frac{5}{2}}\ell^{1+
\frac{1}{k}}t^{\frac{3}{2}-\frac{2}{k}}} \sum_{0<|u| \ll
 U_0 X^\epsilon} (u, s_2^2 s_4
t)^{\frac{1}{2}}
\\& \ll \frac{X^{\frac{1}{2} +
\frac{1}{k}}\|\Psi_{M,0}(\psi_j)\|_1}{s_1^{2+
\frac{1}{k}}s_2^{4 -
\frac{2}{k}}s_3^2s_4^{\frac{5}{2}}\ell^{1+
\frac{1}{k}}t^{\frac{3}{2}-\frac{2}{k}}} \sum_{d | s_2^2 s_4
t} d^{\frac{1}{2}} \sum_{0
< |u| \ll \frac{ U_0 X^\epsilon}{d}} 1.
\end{align*}
For the $L^1$ norm $\|\Psi_{M,0}(\psi_j)\|_1$ recall  (\ref{L_1_bound})
\begin{equation*}
 \|\Psi_{M,0}(\psi_j)\|_1 \ll X^\epsilon M^{-\frac{k-2}{4}}.
\end{equation*}
  Substituting this bound, and
the bound $U_0 \ll \frac{s X^{\frac{k-3}{2}}}{Y_j^{k-1} (s_1
\ell)^{k-3}s_2 t}$ in
(\ref{f_0_u_length}), we obtain
\[
 E_{0,j,e,s,\ell,t} \ll
 \frac{X^{\frac{k}{4} +
\epsilon}}{T^{\frac{k}{2}}s_1^{\frac{k}{2} -
\frac{1}{2}}s_2^3 s_3^2s_4^{\frac{5}{2}}\ell^{1 +
\frac{1}{k}}t^{\frac{3}{2}}}
\]
and thus
\[
 E_{0,j,e} = \sum_{s, \ell, t} E_{0,j,e,s,\ell,t} \ll
\frac{X^{\frac{k}{4} +
\epsilon}}{T^{\frac{k}{2}}}.
\]
Since there are $O(\log X)$ components $\psi_j$ in the
partition of unity, we
deduce that the total contribution of terms $E_{0,j,e}$ to
$\cM$ is
$O\left(X^{\frac{k}{4} + \epsilon}T^{-\frac{k}{2}}\right)$, with an analogous
contribution from the odd component. Combined
with the evaluation
of the diagonal in the previous section, this proves
Proposition
\ref{main_term_prop} in the case $f = 0$.

\subsubsection{Bound for $\Delta_f$, $f \neq 0$}
Following our convention, write $\Delta_{f, j, e}$ to indicate the even
term from $\cM_j$. Bound
\[
 \left|S_{k-2}\left(0, \tilde{f}s_2 s_4; s_2^2 s_4 t\right)\right| \ll
(s_2^2
s_4) (f, t)^{\frac{1}{2}} t^{\frac{1}{2} + \epsilon}
\]
to obtain
\begin{align*}
 \Delta_{f,j,e,s, \ell, t} \ll \frac{X^{\frac{1}{2} +
\frac{1}{k}+\epsilon}}{ s_1^{2+
\frac{1}{k}}s_2^{2 -
\frac{2}{k}}s_3^2s_4^{\frac{3}{2}}\ell^{1+
\frac{1}{k}}t^{\frac{3}{2}-\frac{2}{k}-\epsilon}} (f, t)^{\frac{1}{2}}
\left|\Psi_{M,F}^{1,2}(0,0)\right|.
\end{align*}
Bound
\[
 \left|\Psi_{M,F}^{1,2}(0,0)\right| \leq \|\Psi_{M,0}\|_{1} \ll X^\epsilon
M^{-\frac{k-2}{4}}
\]
to obtain
\begin{align*}
 \Delta_{f,j,e} &\ll X^{\frac{3}{2} -
\frac{k}{4}+\epsilon} Y_j^{\frac{k}{2}-1}\sum_{s=s_1s_2s_3s_4}
\sum_{\substack{(s_1\ell)^{\frac{k-1}{2}}s_2t
\\ \ll
X^{\frac{k-2}{4}+\epsilon}Y_j^{-\frac{k}{2}}}}
\frac{ s_1^{\frac{k}{2} -
\frac{7}{2}}\ell^{\frac{k}{2} -
\frac{5}{2}}(f,t)^{\frac{1}{2}}}{s^2 t^{\frac{1}{2}}}
\\& \ll X^{\frac{3}{2} -
\frac{k}{4}} Y_j^{\frac{k}{2}-1}
\left(X^{\frac{k-2}{4}}Y_j^{-\frac{k}{2}}
\right)^{\max(\frac{1}{2}, \frac{k-3}{k-1})}.
\end{align*}
For $k = 3$ this gives a bound of
\[
 \Delta_f \ll X^{\frac{7}{8} + \epsilon}T^{-\frac{1}{4}}.
\]
For $k \geq 5$ this gives a bound of
\[
 \Delta_f \ll X^{\frac{1}{2} + \frac{1}{2(k-1)}+\epsilon}
Y_j^{\frac{1}{k-1}}
\ll X^{\frac{1}{2} + \frac{1}{k} + \epsilon}.
\]

\subsubsection{Bound for $E_f$, $f \neq 0$}

When $u \neq 0$, bound
\[\left|S_{k-2}\left(\ell s_3 u, \tilde{f}s_2 s_4; s_2^2s_4t\right)
\right|
\ll (u, s_2^2s_4t)^{\frac{1}{2}}(s_2^2 s_4 t)^{\frac{1}{2}+\epsilon}\]
and apply the bound (\ref{FT_F_bound}) of Lemma \ref{FT_formula}
to bound
$\Psi^{1,2}_{M,F}$ by
\[
 \left|\Psi^{1,2}_{M, F}\left(\cdot,
 0\right) \right| \ll \frac{M^{\frac{k}{4} +\frac{1}{2}}}{|F|}
\|\Psi_{M,0}(\psi_j)\|_1.
\]
In view of the bound for the $L^1$ norm
(\ref{L_1_bound}), we have
\[
  \left|\Psi_{M,F}^{1,2}(\cdot, 0)\right| \ll X^\epsilon
\frac{M}{|F|} \ll
\frac{X^{\frac{1}{2}-\frac{1}{k} + \epsilon}}{|f|(s_1 \ell)^{1 - \frac{1}{k}}
(s_2 t)^{\frac{2}{k}}T^2}.
\]
 This obtains
\begin{align*}
E_{f, j, e, s, \ell, t} &\ll \frac{X^{1 + \epsilon}}{|f|T^2
s_1^3s_2^4s_3^2s_4^{\frac{5}{2}}\ell^2
t^{\frac{3}{2}}}\sum_{d | s_2^2 s_4 t} d^{\frac{1}{2}} \sum_{u \ll
\frac{1}{d} \frac{|f| s
X^{\frac{k}{4}-1 + \epsilon}}{T^{\frac{k}{2} -2}}}1 \\&\ll
\frac{
X^{\frac{k}{4} + \epsilon}}{T^{\frac{k}{2}} s_1^2 s_2^3 s_3
s_4^{\frac{3}{2}} \ell^2
t^{\frac{3}{2}}},
\end{align*}
so that, summing over $s, \ell, t$, the
contribution of these terms to
$E_{f,j,e}$ is bounded by $\frac{X^{\frac{k}{4} +
\epsilon}}{T^{\frac{k}{2}}}$.

Combined with the estimate for $\Delta_f$ above and corresponding estimates
in the odd case we obtain  Proposition \ref{main_term_prop} in
the
case $f \neq 0$.

\section{The sieving error term}\label{sieve_section}
The goal of this section is to prove the bound for the sieving error term
claimed in Proposition \ref{error_term_prop}.
The crucial ingredient in the sieve is the following lemma, which associates
to non-square-free $d = d_1 q^2$ and parameterization equation $\ell m^k =
\ell^2
n^2 + t^2 d$, a genuine primitive ideal in the ring of integers of
$\bQ(\sqrt{-d_1})$, and of class lying in a prescribed coset of the $k$-part of
the class group $H(-4d_1)$, with the number of such cosets appearing bounded
by a divisor function of $q$.

Given $q \geq 1$ indicate by
\[
 \sq\left(p_1^{e_1}\cdots
p_r^{e_r}\right) =
p_1^{\left\lfloor \frac{e_1}{2} \right\rfloor}\cdots p_r^{\left\lfloor
\frac{e_r}{2}\right\rfloor}, \qquad
\kr\left(p_1^{e_1}\cdots
p_r^{e_r}\right) =
p_1^{\left\lceil \frac{e_1}{k}\right\rceil}\cdots p_r^{\left\lceil
\frac{e_r}{k}\right\rceil}
\]
the largest number whose square divides $q$, resp. the least $k$th power
divisible by $q$.

\begin{lemma}\label{sieving_lemma}
Let $\left(\ell, m,n,t,q,d\right)\in \left(\zed^+\right)^6$ satisfy $\ell m^k -
\ell^2 n^2 =
t^2 q^2 d$ with $q^2d\equiv 2 \bmod 4$, $d$ square-free, $\left(\ell mn,
t\right)= (\ell, m) = 1$ and $\ell$ square-free.  Set \[ q_1 =
\sq\left(\gcd\left(m^k, n^2\right)\right), \qquad
q_2 = \frac{q}{q_1}.\]  Further, set
also \[q_{10} = \kr(q_1^2).\]  Then define
\[ m' = \frac{m}{q_{10}},
\qquad n' =
  \frac{n}{q_1}, \qquad q' =
\frac{q_{10}^k}{q_{1}^2}.\] The congruence conditions
  $\left(m', \ell  n'\right) = \left(m' n'q', q_2\right) = \left(\ell,
q\right)=1$ hold.  Also, the ideal $\left(q'\right)$ factors
  in $\bQ\left(\sqrt{-d}\right)$ as $\left(q'\right) =
  \mathfrak{q}\overline{\mathfrak{q}}$.  Moreover, there is a
  primitive ideal $\mathfrak{a}$ of $\bQ\left(\sqrt{-d}\right)$ of norm $\ell
  m'$ and solving
  $\mathfrak{q}\mathfrak{a}^k = \ell^{\frac{k-1}{2}} \left(\ell n' + tq_2
\sqrt{-d}\right)$.
\end{lemma}
\begin{proof}
Dividing both sides by $q_1^2$, the equation $\ell m^k - \ell^2 n^2 = t^2q^2d$
may be
rewritten as
\begin{equation}\label{basic}\ell {m'}^k q' - {\ell }^2{n'}^2 = t^2 q_2^2
d.\end{equation}

The condition $(\ell, q) = 1$ follows from $(\ell, m) = 1$ and $\ell$
square-free.  Notice
$\left({m'}^kq', {n'}^2\right)$ is
square-free, and therefore $\left(m', n'\right) = 1$ and also $\left(q',
n'\right)$ is
square-free.  Then $\left(m', \ell  n'\right) = 1$ implies $\left(m',
q_2\right)
= 1$ and
$\left(q', \ell \right) = 1$ implies $\left(q', q_2\right) = 1$ since a
common factor would divide $n'$, but any prime
factor of
$\left(q', n'\right)$ divides $\ell  {m'}^k q'$ only once.  It thus follows that
$\left(n', q_2\right) = 1$, so we have proven all of the congruence conditions.

Equation (\ref{basic}) gives a factorization of ideals
\[\left(\ell  {m'}^k q'\right) = \left(\ell  n' + tq_2
\sqrt{-d}\right)\left(\ell  n' - tq_2
\sqrt{-d}\right)\] in $\bQ\left(\sqrt{-d}\right)$.  Notice $p|\ell  \Rightarrow
p \|
t^2q_2^2d$ so $p|d$, and therefore $(\ell)  = \mathfrak{h}^2$ for some
$\mathfrak{h}$ dividing the different $\mathfrak{d}$.
We claim that $p|q'$ implies $p$ is ramified or split in
$\bQ\left(\sqrt{-d}\right)$.  Indeed, if $p$ is inert then \[\left(p\right) |
\left(\left(\ell  n' +
  tq_2 \sqrt{-d}\right), \left(\ell  n' - tq_2\sqrt{-d}\right)\right) \qquad
\Rightarrow p|n'\] since $p \nmid 2\ell $.  But then $\left(p\right)^2
|{m'}^kq'$, which contradicts $\left({m'}^kq', {n'}^2\right)$ square-free.
Since
all primes dividing $\left(q'\right)$ are ramified or split, we obtain the
factorization $\left(q'\right) = \mathfrak{q}\overline{\mathfrak{q}}$ with
$\mathfrak{q}|\left(\ell  n' + tq_2 \sqrt{-d}\right)$.

Set $\mathfrak{b} = \left(\ell n' + tq_2\sqrt{-d}\right)
\mathfrak{h}^{-1}\mathfrak{q}^{-1}$ so that
$\mathfrak{b}\overline{\mathfrak{b}} = \left(m'\right)^k$.  Note that
\[\left(\mathfrak{b}, \overline{\mathfrak{b}}\right)|\left(2\ell  n', m'\right)
=
\left(1\right)\] and
therefore $\mathfrak{b}$ is primitive, and co-prime to $\mathfrak{d}$.
Therefore there exists primitive ideal $\mathfrak{c}$ satisfying
$\mathfrak{c}^k = \mathfrak{b}$, and furthermore, $\mathfrak{a} =
\mathfrak{h}\mathfrak{c}$ remains primitive.  Clearly
$N\left(\mathfrak{a}\right)
= \ell  m'$ and $\mathfrak{q} \mathfrak{a}^k = \ell^{\frac{k-1}{2}} \left(\ell
n' + tq_2
\sqrt{-d}\right)$ as wanted.
\end{proof}

Before turning to the sieve upper bound, we record bounds regarding
the average number of $k$-torsion elements in the class group.
\begin{proposition}\label{average_k_torsion}
 We have the bounds
\[
 \sum_{\substack{\frac{X}{2} < d < X\\ d \equiv 2\bmod 4\\ \square\text{-free}}}
\sum_{\substack{[(1)] \neq [\fa] \in H(-4d)\\ [\fa]^k = [(1)]}} 1 \ll \left\{
\begin{array}{ccc} X && k = 3\\ X^{\frac{5}{4}} && k = 5\\ X^{\frac{3}{2}
} && k \geq 7 \end{array}\right..
\]
\end{proposition}
\begin{proof} For $k = 3$ this follows from the Davenport-Heilbronn theorem.
For $k = 5$ this follows from the method of Soundararajan \cite{sound}.  When $k
\geq 7$ this is the result of bounding the number of $k$-torsion elements by the
size of the full class group.
\end{proof}

We now prove our basic estimate for the sieve.
\begin{proposition}\label{sieving_proposition}
Let $\phi$ and $\psi$ be non-negative smooth functions having compact support on
$\bR^+$. Let $1 \leq T \ll X^{\frac{1}{2} - \frac{1}{k}}$.  We have the bound
\[\sum_{q > Z} \sum_{\substack{ q^2d \equiv
      2\bmod 4\\ d \; \square\text{-free}}} \sum_{\substack{\left(\ell,
m,n,t\right)\in
\left(\zed^+\right)^4\\
    \left(\ell m n, t\right) = (\ell, t) = 1\\ \ell m^k - \ell^2n^2 = t^2q^2d\\
\ell\;
    \square\text{-free}}}\phi\left(\frac{q^2 d}{X}\right)
\psi\left(\frac{q^2 d}{T^2 \ell^2 m^2  }\right) \ll
\frac{X^{1 + \epsilon}}{TZ} + \frac{X^{\frac{k}{4} +
\epsilon}}{T^{\frac{k}{2}}}.\]
\end{proposition}

\begin{proof} Keep the meaning of $q_{ij}$ etc from Lemma \ref{sieving_lemma},
in particular $q = q_1q_2$ and $q'q_1^2 = q_{10}^k$.
The sum in question is
\begin{align*} &\sum_{\substack{ q = q_{1}q_{2} > Z \\q \text{ odd}\\
  q' q_1^2 = q_{10}^k}}
{\sum_{\substack{ d \equiv 2 \bmod 4\\ \square\text{-free}}}}
\sum_{\substack{(\ell , m', n', t)\in \left(\zed^+\right)^4\\
\cC_6}}\phi\left(\frac{q^2 d}{X}\right)\psi\left(
\frac{q^2 d}{T^2 {\ell }^2 {m'}^2 q_{10}^2  }\right);\\
&\cC_6 = \left\{\begin{array}{l}\ell  \;
\text{square-free},\\ m'
    \text{ odd}\\ (\ell m'n'q', tq_2) = (\ell ,m'q) =
(m',n')=1\\
    \ell  {m'}^k q' - {\ell }^2{n'}^2 = t^2q_2^2
d  \end{array}\right..\end{align*}

Case 1: $\sqrt{d} \gg \ell  m' \asymp \frac{\sqrt{X}}{Tq_{10}}$.

In the first case, set $t' = tq_2$ to obtain, for a suitable non-negative
$\psi_0 \in C_c^\infty(\bR^+)$,
\begin{align*}&\ll X^\epsilon \sum_{\substack{q = q_{1}q_{2}\\
    q\text{ odd}\\ q'q_1^2 = q_{10}^k }}
\sum_{\substack{(\ell ,m',n', t') \in \left(\zed^+ \right)^4, q_2|t'\\ \cC_7}}
  \psi_0\left(\frac{\sqrt{X}}{T {\ell } {m'} q_{10}
}\right);\\
&\cC_7 = \left\{\begin{array}{l}(\ell
m'n'q', t')=
    (m',\ell n')=1\\ \ell {m'}^kq' -{\ell }^2{n'}^2 = {t'}^2d\\ d
    \equiv 2 \bmod 4, \; \text{square-free}\\
    \frac{X}{T^2 q_{10}^2} \ll d <
      \min\left(\frac{X}{Z^2}, \frac{X}{q^2}\right)
\end{array}\right..\end{align*}
Controlling the size of $d$ with a partition of unity, the inner sum is bounded
by (we write $t$ for $t'$)
\begin{align}\label{error_sum} &\sum_{\substack{\max\left(1,
    \frac{X}{T^2q_{10}^2}\right)\\ < e^a <
    \min\left(\frac{X}{Z^2}, \frac{X}{q^2}\right)}}
\sum_{\substack{(\ell q',
      t)=1\\ \ell \; \square\text{-free}\\ (\ell t)^2 \ll
\frac{q'
       X^{\frac{k}{2}}}{e^a T^k q_{10}^{k}}}}
  \sum_{\substack{\left(m', 2\ell t\right)=1\\ m' \asymp
\frac{\sqrt{X}}{Tq_{10} \ell }\\ m' \gg \left(\frac{t^2
            e^a}{\ell  q'}\right)^{\frac{1}{k}}}}\psi_0\left(\frac{\sqrt{X}}{T
{\ell } {m'} q_{10}
}\right)
   \\ \notag& \times \sum_{\substack{\left(n', m't\right)=1\\ \ell {m'}^kq' -
{\ell }^2{n'}^2
       \\ \equiv 2{t}^2 \bmod 4{t}^2}} \sigma^\times\left(\frac{\ell {m'}^kq' -
      {\ell }^2{n'}^2}{t^2 e^a}\right).
\end{align}
Splitting the sum over $n'$ into blocks of length ${t}^2$, this sum
is
\[\ll X^\epsilon\left(O\left(1\right) + \frac{1}{{t}^2}\frac{{t}^2
e^a}{{\ell }^{\frac{3}{2}}
    {m'}^{\frac{k}{2}} {q'}^{\frac{1}{2}}}\right) \ll
X^{\epsilon}\left(O\left(1\right) +
\frac{e^a}{{\ell }^{\frac{3}{2}}
{m'}^{\frac{k}{2}}{q'}^{\frac{1}{2}}} \right).\]
Bounding the sums over $m'$ and $\ell t$ by their length (recall that $q_2|t$),
the
$O\left(1\right)$ term contributes
\begin{align*}&\ll \frac{X^{\frac{k}{4} + \frac{1}{2} +
\epsilon}}{T^{\frac{k}{2}+1}}
\sum_{\substack{Z<q = q_1 q_2  \ll X^{\frac{1}{2}}\\ q' q_1^2 = q_{10}^k}}
\frac{{q'}^{\frac{1}{2}}}{q_{10}^{\frac{k+2}{2}} q_2}
 \sum_{\max\left(1,
    \frac{X}{T^2q_{10}^2}\right) < e^a} \frac{1}{e^{\frac{a}{2}}} \ll
\frac{X^{\frac{k}{4} + \epsilon}}{T^{\frac{k}{2}}}.
\end{align*}
The second term contributes
\begin{align*}
&\ll X^\epsilon\sum_{\substack{Z < q_1q_2 < X^{\frac{1}{2}}\\ q'q_1^2 =
q_{10}^k }}
\frac{1}{{q'}^{\frac{1}{2}}} \sum_{e^a < \frac{X}{Z^2}} e^a
\sum_{\substack{\left(\ell t\right)^2 \ll
\frac{q'
       X^{\frac{k}{2}}}{e^a T^k q_{10}^{k}}\\ q_2 |t}}
\frac{1}{{\ell }^{\frac{3}{2}}}  \sum_{m' \gg
\left(\frac{{t}^2 e^a}{\ell
      q'}\right)^{\frac{1}{k}}} \frac{1}{{m'}^{\frac{k}{2}}}\\
&\ll X^\epsilon\sum_{\substack{Z < q_1q_2 < X^{\frac{1}{2}}\\ q'q_1^2 =
q_{10}^k }}
\frac{1}{{q'}^{\frac{1}{k}}} \sum_{e^a < \frac{X}{Z^2}} e^{(\frac{1}{2} +
\frac{1}{k})a}
\sum_{\substack{\left(\ell t\right)^2 \ll
\frac{q'
       X^{\frac{k}{2}}}{e^a T^k q_{10}^{k}}\\ q_2|t}}
\frac{1}{{\ell }^{1+\frac{1}{k}} {t}^{1 - \frac{2}{k}}}
\\& \ll \frac{X^{\frac{1}{2} + \epsilon}}{T} \sum_{q <
X^{\frac{1}{2}}} \frac{1}{q} \sum_{e^a < \frac{X}{Z^2}} e^{\frac{a}{2}} \ll
\frac{X^{1 + \epsilon}}{TZ}.
\end{align*}

Case 2: $\sqrt{d} \ll \ell  m' \asymp \frac{\sqrt{X}}{Tq_{10}}$.

 Recall from Lemma \ref{norm_bound} that the number of ideals of
$\bQ(\sqrt{-d})$ of a
fixed class, and with norm bounded by $Y \sqrt{d}$ is $\ll(1 + Y)$.  Using
this, we find that  the second case gives
\begin{align*}X^\epsilon\sum_{\substack{q = q_1q_2  > Z\\ q'q_1^2 = q_{10}^k}}
{\sum_{\substack{d
\equiv 2
      \bmod 4 \\ \square\text{-free}\\ d \asymp \frac{X}{q^2}}}}
\sum_{\substack{\left(q'\right) = \mathfrak{q}\overline{\mathfrak{q}}\\ \text{
in
  }\bQ\left(\sqrt{-d}\right)}}
\sum_{\substack{ \mathfrak{a} \text{ primitive in }  \bQ\left(\sqrt{-d}\right)\\
[ \mathfrak{q}\mathfrak{a}^k] = [\left(1\right)],  N\fa \gg\sqrt{d}}}
\psi_0\left(\frac{\sqrt{X}}{T q_{10} N\fa}\right)
\end{align*}
The support of $\psi_0$ imposes $q_{10} \ll \frac{\sqrt{X}}{T}$.
Also, knowing the ideal $\fa$ we recover
$q_2t$, and hence $q_2$ up to a divisor function.  Putting these together, we
obtain
\begin{align*}
\ll \frac{X^{\frac{1}{2} + \epsilon}}{T} \sum_{\substack{ q_{10} \ll
\frac{\sqrt{X}}{T}}}\frac{1}{q_{10}}
\sum_{\substack{d \equiv 2 \bmod 4, \square\text{-free}\\ d < \min(
\frac{X}{Z^2}, \frac{X}{T^2 q_{10}^2})}} \frac{1}{\sqrt{d}}
\sum_{\substack{[\fa] \in H(-4d)\\ [\fa]^k = [(1)]}} 1
\end{align*}
Substituting the bounds for the average number of $k$-torsion elements
(Proposition \ref{average_k_torsion}) we obtain a bound of $\ll \frac{X^{1
+ \epsilon}}{TZ}$ for $k = 3$, $\ll \frac{X^{\frac{5}{4} +
\epsilon}}{T^{\frac{5}{2}}}$ for $k = 5$,  $\ll \frac{X^{\frac{3}{2} +
\epsilon}}{T^3}$ for $k \geq 7$.  This completes the proof.
\end{proof}

\begin{proof}[Proof of Proposition \ref{error_term_prop}]
 We bound $\mathscr{E}_j$ by
\begin{align*} \mathscr{E}_j &\leq \sum_{\substack{\ell,m,t \in \zed^+, n \in
\zed\\ (\ell mn,t) = (\ell, m)
= 1\\ \ell\; \square\text{-free}\\ \ell m^k  - \ell^2n^2 \equiv 2 t^2
\bmod
4t^2}} \Phi\left(\ell m,
\ell^{\frac{k+1}{2}}n, \ell^{\frac{k-1}{2}}t\Big| \psi_j \right)
\sum_{\substack{s^2|\frac{\ell
m^k-\ell^2n^2}{t^2}, s > Z}} 1
\\
& \ll X^\epsilon\sum_{q > Z} \sum_{\substack{ q^2d \equiv
      2\bmod 4\\ d\, \square\text{-free}}} \sum_{\substack{\left(\ell,
m,n,t\right)\in
\left(\zed^+\right)^4\\
    \left(\ell m n, t\right) = 1, \ell\;
    \square\text{-free}\\ \ell m^k - \ell^2n^2 =
t^2q^2d}}\phi\left(\frac{q^2
d}{X}\right)
\psi_j\left(\frac{\ell m^k - \ell^2 n^2}{Y_j^2 \ell^2 m^2 t^2 }\right),
\end{align*}
which reduces to the sum estimated in Proposition \ref{sieving_proposition}.
\end{proof}

\subsection*{Acknowledgments}
The problem of proving equidistribution was suggested to me by graduate advisor
Soundararajan.  The possibility of calculating a negative secondary term related
to
$k$-torsion was suggested to me by Akshay Venkatesh.

\end{document}